\documentclass[preprint,12pt]{elsarticle}

\usepackage[top=1in, bottom=1in, left=1in, right=1in]{geometry}
\usepackage{amsmath,amssymb,amsfonts,amsthm,bm,dsfont,graphicx,subfigure,epstopdf,hyperref,color,float,multirow}

\numberwithin{equation}{section}
\newtheorem{theorem}{Theorem}[section]

\newtheorem{remark}{Remark}[section]

\newcommand{\rd}{{\rm d}}
\graphicspath{{FIGS/}}

\begin{document}

\begin{frontmatter}
\title{ {\Large
Partitioned AVF methods}}
\author{Wenjun Cai$^{a}$, Haochen Li$^{b}$, Yushun Wang$^{a,*}$}
\address{
$^{a}$ Jiangsu Key Laboratory for NSLSCS, School of Mathematical Sciences, Nanjing Normal University, Nanjing, $210023$, China\\
$^{b}$ School of Mathematical Sciences, Peking University, Beijing, 100871, China\\}

\begin{abstract}
The classic second-order average vector field (AVF) method can exactly preserve the energy for Hamiltonian ordinary differential equations and partial differential equations. However, the AVF method inevitably leads to fully-implicit nonlinear algebraic equations for general nonlinear systems. To address this drawback and maintain the desired energy-preserving property, a first-order partitioned AVF method is proposed which first divides the variables into groups and then applies the AVF method step by step. In conjunction with its adjoint method we present the partitioned AVF composition method and plus method respectively to improve its accuracy to second order. Concrete schemes for two classic model equations are constructed with semi-implicit, linear-implicit properties that make considerable lower cost than the original AVF method. Furthermore, additional conservative property can be generated besides the conventional energy preservation for specific problems. Numerical verification of these schemes further conforms our results.

\end{abstract}
\end{frontmatter}

\begin{figure}[b]
\small \baselineskip=10pt
\rule[2mm]{1.8cm}{0.2mm} \par
$^{*}$Corresponding author.\\
E-mail address: wangyushun@njnu.edu.cn (Y. Wang).
\end{figure}

\pagestyle{myheadings}
\markboth{\hfil 
   \hfil \hbox{}}
{\hbox{} \hfil 
Wenjun Cai,  Yushun Wang, Haochen Li \hfil}

\section{Introduction}
\label{sec:1}

We consider a Hamiltonian system of the form
\begin{equation}\label{eq:1-1}
	\dot{z}=f(z)=S_{m}\nabla H(z),
\end{equation}
where $z\in \mathbb{R}^{m}$, $S_{m}$ is a $m\times m$ skew-symmetric constant matrix, and $H(z)$ is a Hamiltonian assumed to be sufficiently differentiable. The most relevant characters of the system \eqref{eq:1-1} are the symplecticity and energy conservation along any continuous flows. However, under discrete circumstances, numerical integrations cannot fulfill both properties at the same time for arbitrary Hamiltonian $H(z)$ \cite{gm88,cfm06}. As a consequence, symplectic integrators and energy-preserving schemes constitute two prominent lines to construct geometric numerical integration methods for Hamiltonian systems. 

In contrast to earlier commenced symplectic integrators (see the monographs \cite{sc94,lr04,hlw06}), the energy-preserving methods have only emerged in recent decades. An initial and natural approach is to force a nonconservative scheme to lie on a proper manifold representing a first integral of the original system, leading to the projection and symmetric projection strategy \cite{h00,hw06}. The development of the discrete gradient method brings a significant improvement for the construction of energy-preserving methods, which was perhaps first treated in a systematic way by Gonzalez \cite{g96}, see also McLachlan et al. \cite{mqr99}. For systems with polynomial type Hamiltonian functions, Brugnano et al. proposed Hamiltonian boundary value methods (HBVMs) in a series of papers \cite{bit09,bit10,bit11,bi12} that can ensure energy conservation of polynomial Hamiltonians of any high degree. When considering Hamiltonian PDEs such that the gradient of Hamiltonian is replaced by its variation derivative,  Furihata et al. \cite{f99,fm11} developed the discrete variational derivative method (DVDM)  from a finite difference perspective. Later, this method is extended by Matsuo and referred as the discrete partial derivative method (DPDM) \cite{m08}, which focuses on the Galerkin framework.  Recently, a new class of energy-preserving methods called average vector field (AVF) method was introduced in \cite{qm08}, which can be characterized as a specific discrete gradient method by taking a mean-value discrete gradient as the discrete counterpart of the gradient operator \cite{mqr99}. Also the AVF method is a limit case of HBVMs,  setting the stage of Butcher tableau to infinity \cite{bit12}. Furthermore, such a method is equivalent to the DPDM providing the Hamiltonian does not contain cross-product terms of solution and its partial derivatives \cite{m08}. Additionally, the AVF method has been analysed in the framework of B-series \cite{qm08} and prompted intensive studies on energy-preserving B-series methods \cite{cmmoqw09,h10,cmoq10,cos14}.

For the Hamiltonian system \eqref{eq:1-1}, the AVF method is defined by
\begin{equation}\label{eq:1-2}
\frac{z^{n+1}-z^{n}}{\tau}=S_{m}\int_{0}^{1}\nabla H(\xi z^{n+1}+(1-\xi)z^{n})\rd\xi,
\end{equation}
where $\tau$ denotes the time step. Using the fundamental theorem of calculus and the skew-symmetry of $S_{2m}$, the Hamiltonian energy is precisely conserved at every time step, that is $H(z^{n+1})=H(z^n)$. 
One of remarkable advantages of the AVF method is that it only requires evaluations of the vector field. For polynomial Hamiltonians, the integral can be evaluated exactly, and the implementation is comparable to that of the implicit mid-point rule. When the Hamiltonian energy is a quadratic function, the resulting AVF scheme is linearly implicit and therefore can be efficiently solved. 
But this is not the case to reflect the merit of the AVF method since any symplectic integrator can also achieve the energy conservation for quadratic Hamiltonians \cite{hlw06}. Under most circumstances, the evaluation of the integration in \eqref{eq:1-2} leads to a nonlinear function of $z^{n+1}$ which further constitutes a fully implicit numerical scheme. The iterative processes are then inevitably required such that the computational complexity will be evidently increased, especially for the application on Hamiltonian PDEs.

The main aim of this paper is to construct a more efficient AVF based method. Instead of imposing the mean-value discrete gradient straightly along the path jointing $z^n$ and $z^{n+1}$, we first divide a single path into several subpaths by grouping the components of $z$ and then apply the  mean-value discrete gradient one group at a time. We denote this method as a partitioned AVF (PAVF) method which can also automatically preserve arbitrary Hamiltonian energy of system \eqref{eq:1-1}. Comparing with the AVF method \eqref{eq:1-2}, the resulting schemes of the partitioned AVF method are much simpler according to the concrete expressions which reduces the original fully implicit schemes to semi-implicit or linearly implicit schemes and therefore significantly improve the computational efficiency. Although such a method is only of first-order accuracy, in conjunction with its adjoint we further present the partitioned AVF composition (PAVF-C) method and the partitioned AVF plus (PAVF-P) method that both achieve second-order accuracy and energy-preserving property. Thanks to the great advantage of the PAVF method, even by its composition the PAVF-C method still has much lower cost than the direct AVF method. Although the computational efficiency of the PAVF-P method is comparable to the AVF method, it may possess additional conservative quantities that the AVF method cannot preserve.

This paper is organized as follows. In section 2, we present the partitioned AVF method and derive its energy-preserving property. The adjoint of this method as well as the induced partitioned AVF composition and plus methods are also proposed. Numerical schemes for Hamiltonian ODEs and PDEs are constructed and tested in section 3, including the H\'{e}non-Heiles system and the Klein-Gordon-Schr\"{o}dinger equation. We compare the original AVF method and our partitioned version methods in both the invariant preservation and computational cost. Concluding remarks are given in Section 4.

\section{The partitioned AVF methods}
\label{sec:2}

For illustration, we consider the Hamiltonian system \eqref{eq:1-1} when $m$ is an even number, denoting $m=2d$. Without generality the grouping strategy is simply choosing  in sequential order, i.e., $z=(p,q)^T=(z_{1},z_{2},...,z_{d} ;z_{d+1},z_{d+2},...,z_{m})^{T}$. Accordingly, the original system \eqref{eq:1-1} can be rewritten as
\begin{align}\label{eq:2-4}
\left(
       \begin{array}{c}
              \dot{p}   \\
              \dot{q}
             \end{array}
\right)
=S_{2d}
\left(
       \begin{array}{c}
              H_{p}(p,q)   \\
              H_{q}(p,q)
             \end{array}
\right),\quad
p,q\in \mathbb{R}^{d}.
\end{align}
The present Hamiltonian $H(p,q)$ is still conserved along any continuous flow, that is 
\begin{equation*}
\frac{\rd H(p(t),q(t))}{\rd t}=H_{p}(p,q)^{T}\dot{p}+H_{q}(p,q)^{T}\dot{q}=\nabla H(p,q)^{T}S_{2d}\nabla H(p,q)=0.
\end{equation*}
Then the so-called partitioned AVF (PAVF) method for the Hamiltonian system \eqref{eq:2-4} is defined by
\begin{align}\label{eq:2-5}
\frac{1}{\tau}\left(
       \begin{array}{c}
              p^{n+1}-p^{n} \\
              q^{n+1}-q^{n}
             \end{array}
\right)
=S_{2d}
\left(
       \begin{array}{c}
              \int_{0}^{1}H_{p}(\xi p^{n+1}+(1-\xi)p^{n},q^{n})\rd\xi \\
              \int_{0}^{1}H_{q}(p^{n+1},\xi q^{n+1}+(1-\xi)q^{n})\rd\xi
             \end{array}
\right).
\end{align}

\begin{theorem}\label{thm:2-1}
The PAVF method \eqref{eq:2-5} preserves the Hamiltonian $H(p,q)$ of the system \eqref{eq:2-4} exactly, and satisfies
\begin{equation*}
\frac{1}{\tau}(H(z^{n+1})-H(z^{n}))=0.
\end{equation*}
\end{theorem}
\begin{proof}
Taking the scalar product with $(\int_{0}^{1}H_{p}(\xi p^{n+1}+(1-\xi)p^{n},q^{n})^{T}\rd\xi,\int_{0}^{1}H_{q}(p^{n+1},\xi q^{n+1}+(1-\xi)q^{n})^{T}\rd\xi)^T$ on both sides of \eqref{eq:2-5}, using the Fundamental Theorem of Calculus and the skew-symmetry of $S_{2d}$, we obtain
\begin{align*}
0&=\frac{1}{\tau}\int_{0}^{1}H_{p}(\xi p^{n+1}+(1-\xi)p^{n},q^{n})^{T}\rd\xi(p^{n+1}-p^{n})\\
&\qquad+\frac{1}{\tau}\int_{0}^{1}H_{q}(p^{n+1},\xi q^{n+1}+(1-\xi)q^{n})^{T}\rd\xi(q^{n+1}-q^{n})\\
&=\frac{1}{\tau}\int_{0}^{1}\frac{\rd}{\rd\xi}\big[H(\xi p^{n+1}+(1-\xi)p^{n},q^{n})+H(p^{n+1},\xi q^{n+1}+(1-\xi)q^{n})\big]\rd\xi\\
&=\frac{1}{\tau}\big[H(p^{n+1},q^{n})-H(p^{n},q^{n})+H(p^{n+1},q^{n+1})-H(p^{n+1},q^{n})\big]\\
&=\frac{1}{\tau}(H(p^{n+1},q^{n+1})-H(p^{n},q^{n}))=0.
\end{align*}
\end{proof}

\begin{remark}\label{rmk:2-1}
	The PAVF methods \eqref{eq:2-5} are one-step methods of order one.
\end{remark}


\begin{remark}\label{rmk:2-3}
	If the  system \eqref{eq:2-4} is a separable Hamiltonian system of the form $H(p,q)=H_{1}(p)+H_{2}(q)$,
	then the PAVF method \eqref{eq:2-5} is equivalent to the AVF method \eqref{eq:1-2}. However, the Hamiltonian system \eqref{eq:2-4} is not separable in general which can just reflect the major advantage of the PAVF method. For illustration, consider a specific case when $H_p(p,q)=pq$, then the conventional AVF method will lead to a nonlinear term $\frac{1}{3}p^{n+1}q^{n+1}+\frac{1}{6}(p^nq^{n+1}+p^{n+1}q^n)+\frac{1}{3}p^nq^n$ and therefore requires an iteration process for the associated nonlinear equations. While for the PAVF method, we can divide $p$, $q$ into two individual groups and apply the AVF method only for one group during each step. The corresponding scheme is linearly implicit at least for a part of the entire equation system which dramatically reduces the iteration scale or even avoid the iteration completely. For the general Hamiltonian system, the grouping strategy is to separate the variables appearing as cross-product terms in the gradient of the Hamiltonian and classify them into different groups. As a consequence, the resulting PAVF method is not unique but can all preserve the Hamiltonian with lower computational cost.

\end{remark}

We denote the PAVF method \eqref{eq:2-5} as $\Phi_\tau$, then its adjoint method $\Phi^*_\tau$ can be obtained as follows
\begin{align}\label{eq:2-5m}
\frac{1}{\tau}\left(
\begin{array}{c}
p^{n+1}-p^{n} \\
q^{n+1}-q^{n}
\end{array}
\right)
=S_{2d}
\left(
\begin{array}{c}
\int_{0}^{1}H_{p}(\xi p^{n+1}+(1-\xi)p^{n},q^{n+1})\rd\xi \\
\int_{0}^{1}H_{q}(p^{n},\xi q^{n+1}+(1-\xi)q^{n})\rd\xi
\end{array}
\right).
\end{align}
Obviously, the adjoint method $\Phi_{\tau}^{*}$ \eqref{eq:2-5m} belongs to the PAVF method, only with a reversed path order to that of $\Phi_\tau$. Consequently, the adjoint method \eqref{eq:2-5m} can also exactly preserve the Hamiltonian $H(z)$ and possess all the properties of the PAVF method.

In conjunction with the adjoint PAVF method, we can define the following PAVF composition (PAVF-C) method
\begin{equation}\label{eq:2-9}
\Psi_{\tau}:=\Phi_{\frac{\tau}{2}}^{*}\circ\Phi_{\frac{\tau}{2}},
\end{equation}
and partitioned AVF plus (PAVF-P) method
\begin{equation}\label{eq:2-10}
\hat{\Psi}_{\tau}:=\frac{1}{2}(\Phi_{\tau}^{*}+\Phi_{\tau}),
\end{equation}
respectively.

\begin{theorem}\label{thm:2-3}
	The PAVF-C method $\Psi_{\tau}$ \eqref{eq:2-9}  and PAVF-P method $\hat{\Psi}_{\tau}$ \eqref{eq:2-10}  are both second-order methods which can conserve the energy of the Hamiltonian system \eqref{eq:2-6} exactly.
\end{theorem}

\begin{proof}
	It is clear that these two methods are both symmetric, the second-order accuracy is thereby obviously. Using the fact that either the PAVF method or its adjoint can exactly preserve the same Hamiltonian, the operators of their composition \eqref{eq:2-9} and plus \eqref{eq:2-10} still inherit this property that ends the proof. 
\end{proof}

Next, we consider a more general case of the Hamiltonian system \eqref{eq:1-1} with $m=\tilde{m}d$ which very often associates to the semi-discrete system of Hamiltonian PDEs. The corresponding grouping strategy is the same as that in Remark \ref{rmk:2-3} and for convenience we just take $z=(\tilde{z}_1,\tilde{z}_2,...,\tilde{z}_{\tilde{m}-1},\tilde{z}_{\tilde{m}})$ where $\tilde{z}_k\in\mathbb{R}^d$, $k=1,2,...,\tilde{m}$. Then we can reform the Hamiltonian system \eqref{eq:1-1} as
\begin{align}\label{eq:2-6}
\left(
\begin{array}{c}
\dot{\tilde{z}}_1   \\
\dot{\tilde{z}}_2\\
\vdots\\
\dot{\tilde{z}}_{\tilde{m}-1}   \\
\dot{\tilde{z}}_{\tilde{m}}
\end{array}
\right)
=S_{\tilde{m}d}
\left(
\begin{array}{c}
H_{\tilde{z}_{1}}(\tilde{z}_{1},\tilde{z}_{2},...,\tilde{z}_{\tilde
	{m}-1},\tilde{z}_{\tilde
	{m}})  \\
H_{\tilde{z}_{2}}(\tilde{z}_{1},\tilde{z}_{2},...,\tilde{z}_{\tilde
	{m}-1},\tilde{z}_{\tilde
	{m}})\\
\vdots\\
H_{\tilde{z}_{{\tilde
			{m}-1}}}(\tilde{z}_{1},\tilde{z}_{2},...,\tilde{z}_{\tilde
	{m}-1},\tilde{z}_{\tilde
	{m}})  \\
H_{\tilde{z}_{{\tilde
			{m}}}}(\tilde{z}_{1},\tilde{z}_{2},...,\tilde{z}_{\tilde
	{m}-1},\tilde{z}_{\tilde
	{m}})  \\
\end{array}
\right).
\end{align}
For this general case, the PAVF method is defined by
\begin{align}\label{eq:2-7}
\frac{1}{\tau}\left(
       \begin{array}{c}
              \tilde{z}_{1}^{n+1}-\tilde{z}_{1}^{n} \\
              \tilde{z}_{2}^{n+1}-\tilde{z}_{2}^{n} \\
              \vdots \\
              \tilde{z}_{\tilde{m}-1}^{n+1}-\tilde{z}_{\tilde{m}-1}^{n} \\
             \tilde{z}_{\tilde{m}}^{n+1}-\tilde{z}_{\tilde{m}}^{n}
             \end{array}
\right)
=S_{\tilde{m}d}
\left(
       \begin{array}{c}
              \int_{0}^{1}H_{\tilde{z}_{1}}(\xi \tilde{z}_{1}^{n+1}+(1-\xi)\tilde{z}_{1}^{n},\tilde{z}_{2}^{n},...,\tilde{z}_{\tilde
              	{m}-1}^{n},\tilde{z}_{\tilde
              {m}}^{n})\rd\xi \\
              \int_{0}^{1}H_{\tilde{z}_{2}}(\tilde{z}_{1}^{n+1},\xi \tilde{z}_{2}^{n+1}+(1-\xi)\tilde{z}_{2}^{n},...,\tilde{z}_{\tilde
              	{m}-1}^{n},\tilde{z}_{\tilde
              	{m}}^{n})\rd\xi \\
              \vdots \\
              \int_{0}^{1}H_{\tilde{z}_{\tilde
              		{m}-1}}(\tilde{z}_{1}^{n+1},\tilde{z}_{2}^{n+1},...,\xi \tilde{z}_{\tilde
              	{m}-1}^{n+1}+(1-\xi)\tilde{z}_{\tilde
              	{m}-1}^{n},z_{\tilde
              	{m}}^{n})\rd\xi \\
              \int_{0}^{1}H_{\tilde{z}_{\tilde
              		{m}}}(\tilde{z}_{1}^{n+1},\tilde{z}_{2}^{n+1},...,\tilde{z}_{\tilde
              	{m}-1}^{n+1},\xi \tilde{z}_{\tilde
              	{m}}^{n+1}+(1-\xi)\tilde{z}_{\tilde
              	{m}}^{n})\rd\xi
             \end{array}
\right).
\end{align}

\begin{theorem}\label{thm:2-2}
The PAVF method \eqref{eq:2-7} preserves the Hamiltonian $H(z)$ of the general system \eqref{eq:2-6} exactly, and satisfies
\begin{equation*}
\frac{1}{\tau}(H(z^{n+1})-H(z^{n}))=0.
\end{equation*}
\end{theorem}
\begin{proof}
In analogy to the proof of above theorem (\ref{thm:2-1}), we can get
\begin{align*}
0=&\frac{1}{\tau}\sum_{k=1}^{\tilde{m}}\big[\int_{0}^{1}H_{\tilde{z}_{k}}(\tilde{z}_{1}^{n+1},...,\tilde{z}_{k-1}^{n+1},\xi \tilde{z}_{k}^{n+1}+(1-\xi)\tilde{z}_{k}^{n},\tilde{z}_{k+1}^{n},...,\tilde{z}_{\tilde{m}}^{n})d\xi\big]^{T}
(\tilde{z}_{k}^{n+1}-\tilde{z}_{k}^{n})\\
=&\frac{1}{\tau}\sum_{k=1}^{\tilde{m}}(H(\tilde{z}_{1}^{n+1},...,\tilde{z}_{k-1}^{n+1},\tilde{z}_{k}^{n+1},\tilde{z}_{k+1}^{n},...,\tilde{z}_{\tilde{m}}^{n})-
(\tilde{z}_{1}^{n+1},...,\tilde{z}_{k-1}^{n+1},\tilde{z}_{k}^{n},\tilde{z}_{k+1}^{n},...,\tilde{z}_{m}^{n}))\\
=&\frac{1}{\tau}(H(\tilde{z}_{1}^{n+1},\tilde{z}_{2}^{n+1},\tilde{z}_{3}^{n+1},...,\tilde{z}_{\tilde{m}}^{n+1})-
(\tilde{z}_{1}^{n},\tilde{z}_{2}^{n},\tilde{z}_{3}^{n},...,\tilde{z}_{\tilde{m}}^{n})).
\end{align*}
\end{proof}

We can also derive the its adjoint scheme as
\begin{align}\label{eq:2-8}
\frac{1}{\tau}\left(
\begin{array}{c}
\tilde{z}_{1}^{n+1}-\tilde{z}_{1}^{n} \\
\tilde{z}_{2}^{n+1}-\tilde{z}_{2}^{n} \\
\vdots \\
\tilde{z}_{\tilde{m}-1}^{n+1}-\tilde{z}_{\tilde{m}-1}^{n} \\
\tilde{z}_{\tilde{m}}^{n+1}-\tilde{z}_{\tilde{m}}^{n}
\end{array}
\right)
=S_{\tilde{m}d}
\left(
\begin{array}{c}
\int_{0}^{1}H_{\tilde{z}_{1}}(\xi \tilde{z}_{1}^{n+1}+(1-\xi)\tilde{z}_{1}^{n},\tilde{z}_{2}^{n+1},...,\tilde{z}_{\tilde
	{m}-1}^{n+1},\tilde{z}_{\tilde
	{m}}^{n+1})\rd\xi \\
\int_{0}^{1}H_{\tilde{z}_{2}}(\tilde{z}_{1}^{n},\xi \tilde{z}_{2}^{n+1}+(1-\xi)\tilde{z}_{2}^{n},...,\tilde{z}_{\tilde
	{m}-1}^{n+1},\tilde{z}_{\tilde
	{m}}^{n+1})\rd\xi \\
\vdots \\
\int_{0}^{1}H_{\tilde{z}_{\tilde
		{m}-1}}(\tilde{z}_{1}^{n},\tilde{z}_{2}^{n},...,\xi \tilde{z}_{\tilde
	{m}-1}^{n+1}+(1-\xi)\tilde{z}_{\tilde
	{m}-1}^{n},z_{\tilde
	{m}}^{n+1})\rd\xi \\
\int_{0}^{1}H_{\tilde{z}_{\tilde
		{m}}}(\tilde{z}_{1}^{n},\tilde{z}_{2}^{n},...,\tilde{z}_{\tilde
	{m}-1}^{n},\xi \tilde{z}_{\tilde
	{m}}^{n+1}+(1-\xi)\tilde{z}_{\tilde
	{m}}^{n})\rd\xi
\end{array}
\right),
\end{align}
Thereafter, the PAVF-C and PAVF-P methods for general Hamiltonian system can be obtained.

Notice that the discrete gradient method \cite{g96,mqr99} is a special case of the proposed PAVF method, which generally has the form
\begin{equation}\label{eq:2-1}
\frac{z^{n+1}-z^{n}}{\tau}=S_{m} \overline{\nabla} H(z^{n+1},z^{n})
\end{equation}
for the system \eqref{eq:1-1}. However, there are various definitions of the discrete gradient $\overline{\nabla}H(z^{n+1},z^n)$ \cite{mqr99}. One can choose the coordinate increment discrete gradient \cite{ia88} as
\begin{equation}\label{eq:2-2}
\overline{\nabla} H(\bar{z},z):=
\left(
\begin{array}{c}
\dfrac{H(\bar{z}_{1},z_{2},z_{2},...,z_{m})-H(z_{1},z_{2},z_{3},...,z_{m})}{\bar{z}_{1}-z_{1}}   \\
\dfrac{H(\bar{z}_{1},\bar{z}_{2},z_{3},...,z_{m})-H(\bar{z}_{1},z_{2},z_{3},...,z_{m})}{\bar{z}_{2}-z_{2}}   \\
\vdots\\
\dfrac{H(\bar{z}_{1},...,\bar{z}_{m-2},\bar{z}_{m-1},z_{m})-H(\bar{z}_{1},...,\bar{z}_{m-2},z_{m-1},z_{m})}{\bar{z}_{m-1}-z_{m-1}}   \\
\dfrac{H(\bar{z}_{1},...,\bar{z}_{m-2},\bar{z}_{m-1},\bar{z}_{m})-H(\bar{z}_{1},...,\bar{z}_{m-2},\bar{z}_{m-1},z_{m})}{\bar{z}_{m}-z_{m}}
\end{array}
\right),
\end{equation}
with the notation $\bar{z}=z^{n+1}$, $z=z^{n}$, which can be further regarded as vector of the means of the tangential components of $\nabla H$ along each of the $m$ segments of the path joining $\bar{z}$ and $z$ by incrementing the coordinates one at a time \cite{mqr99}, that is 
\begin{equation}\label{eq:2-3}
(\overline{\nabla}H(\bar{z},z))_k=\int_{0}^{1}\frac{\partial H}{\partial z_k}(\bar{z}_1,...,\bar{z}_{k-1}, \xi\bar{z}_k+(1-\xi)z_k, z_{k+1},...,z_m)\rd\xi.
\end{equation}
The PAVF method is directly inspired by the discrete gradient method \eqref{eq:2-2} and its equivalent form \eqref{eq:2-3}. The major difference is that instead of dividing the path along $z$ to $\bar{z}$ into $m$ segments, we first group the components according to the cross-product terms in the gradient of the Hamiltonian and then apply the processes of incrementing the group one at a time. As a consequence, the choices of the grouping strategy are flexible, which can lead to more  efficient energy-preserving schemes.

In the next section, we will apply the class of PAVF methods on both ODEs and PDEs, and construct concrete energy-preserving schemes.

\section{Numerical examples}
\label{sec:3}

\subsection{H\'{e}non-Heiles system}
\label{sec:3-1}

Consider the H\'{e}non-Heiles system
\begin{equation}\label{eq:3-1}
\dot{z}=J\nabla H(z),\quad H(z)=\frac{1}{2}(q_{1}^{2}+q_{2}^{2}+p_{1}^{2}+p_{2}^{2})+q_{1}^{2}q_{2}-\frac{1}{3}q_{2}^{3},
\end{equation}
where $z=(q_1,q_2,p_1,p_2)^T$ and $J=\left(
\begin{array}{cc}
0 & I   \\
-I & 0  
\end{array}
\right)$ with $I$ a $2\times 2$ identity matrix. This model was created for describing stellar motion, followed for a very long time, inside the gravitational potential of a galaxy with cylindrical \cite{hlw06}. The H\'{e}non-Heiles system \eqref{eq:3-1} has a finite energy of escape $H_{esc}$ which is equal to $1/6$. For values of energy $H<H_{esc}$, the equipotential curves of the system are close thus making escape impossible. However, for energy greater than $H_{esc}$, the equipotential curves open and three exit channels appear through which the test particles may escape to infinity \cite{zotos}. 

\subsubsection{Derivation of the partitioned AVF schemes}

In contrast, we first present the conventional second-order AVF method \eqref{eq:1-2} for the H\'{e}non-Heiles system which can be derived as follows
\begin{equation}\label{eq:3-2}
\begin{aligned}
&\frac{1}{\tau}(q_{1}^{n+1}-q_{1}^{n})=\frac{1}{2}(p_{1}^{n+1}+p_{1}^{n}), \\
&\frac{1}{\tau}(q_{2}^{n+1}-q_{2}^{n})=\frac{1}{2}(p_{2}^{n+1}+p_{2}^{n}), \\
&\frac{1}{\tau}(p_{1}^{n+1}-p_{1}^{n})=-(\frac{1}{2}(q_{1}^{n+1}+q_{1}^{n})+
\frac{1}{3}(q_{1}^{n+1}q_{2}^{n+1}+4q_{1}^{n+\frac{1}{2}}q_{2}^{n+\frac{1}{2}}+q_{1}^{n}q_{2}^{n})), \\
&\frac{1}{\tau}(p_{2}^{n+1}-p_{2}^{n})=
-\frac{1}{2}(q_{2}^{n+1}+q_{2}^{n})+\frac{1}{3}((q_{2}^{n+1})^{2}+q_{2}^{n+1}q_{2}^{n}+(q_{2}^{n})^{2}-(q_{1}^{n+1})^{2}-q_{1}^{n+1}q_{1}^{n}-(q_{1}^{n})^{2}).
\end{aligned}
\end{equation}
Notice that in this example, the PAVF method for the system \eqref{eq:3-1} is actually equivalent to the discrete gradient method \eqref{eq:2-2} with the corresponding discrete gradient defined by \eqref{eq:2-3}, which is a special case from the proposed general PAVF method \eqref{eq:2-6} when $d=1, m=4$. We can write down the PAVF method for the H\'{e}non-Heiles system \eqref{eq:3-1} as
\begin{equation}\label{eq:3-3}
\begin{aligned}
&\frac{1}{\tau}(q_{1}^{n+1}-q_{1}^{n})=\frac{1}{2}(p_{1}^{n+1}+p_{1}^{n}), \\ 
&\frac{1}{\tau}(q_{2}^{n+1}-q_{2}^{n})=\frac{1}{2}(p_{2}^{n+1}+p_{2}^{n}), \\ 
&\frac{1}{\tau}(p_{1}^{n+1}-p_{1}^{n})=-(\frac{1}{2}(q_{1}^{n+1}+q_{1}^{n})+(q_{1}^{n+1}+q_{1}^{n})q_{2}^{n}), \\
&\frac{1}{\tau}(p_{2}^{n+1}-p_{2}^{n})=
-(\frac{1}{2}(q_{2}^{n+1}+q_{2}^{n})+(q_{1}^{n+1})^{2})+\frac{1}{3}((q_{2}^{n+1})^{2}+q_{2}^{n+1}q_{2}^{n}+(q_{2}^{n})^{2}).
\end{aligned}
\end{equation}
It is clear that the PAVF method \eqref{eq:3-3} is simpler than the original one \eqref{eq:3-2} at first glance. Furthermore, we can find the former scheme is semi-implicit, that is given $(q_1^n,q_2^n,p_1^n,p_2^n)$ the values $q_1^{n+1}$ and $p_1^{n+1}$ can be calculated from the first and third equations of \eqref{eq:3-3} explicitly, leaving only $q_2^{n+1}, p_2^{n+1}$ being solved by numerical iterations. In contrast, the conventional AVF method \eqref{eq:3-2} has to run an iteration involving all the variable $(q_1^{n+1},q_2^{n+1},p_1^{n+1},p_2^{n+1})$ which double the computational scale of the PAVF method \eqref{eq:3-3}. Therefore, the PAVF method is apparently more effective than the AVF method.

In conjunction with the adjoint method of \eqref{eq:3-3} we can derive the PAVF-C method 
\begin{equation}\label{eq:3-4}
\begin{aligned}
&\frac{2}{\tau}(q_{1}^{*}-q_{1}^{n})=\frac{1}{2}(p_{1}^{*}+p_{1}^{n}), \\
&\frac{2}{\tau}(q_{2}^{*}-q_{2}^{n})=\frac{1}{2}(p_{2}^{*}+p_{2}^{n}), \\
&\frac{2}{\tau}(p_{1}^{*}-p_{1}^{n})=-(\frac{1}{2}(q_{1}^{*}+q_{1}^{n})+(q_{1}^{*}+q_{1}^{n})q_{2}^{n}), \\
&\frac{2}{\tau}(p_{2}^{*}-p_{2}^{n})=
-(\frac{1}{2}(q_{2}^{*}+q_{2}^{n})+(q_{1}^{*})^{2})+\frac{1}{3}((q_{2}^{*})^{2}+q_{2}^{*}q_{2}^{n}+(q_{2}^{n})^{2}),\\
&\frac{2}{\tau}(q_{1}^{n+1}-q_{1}^{*})=\frac{1}{2}(p_{1}^{n+1}+p_{1}^{*}), \\
&\frac{2}{\tau}(q_{2}^{n+1}-q_{2}^{*})=\frac{1}{2}(p_{2}^{n+1}+p_{2}^{*}), \\
&\frac{2}{\tau}(p_{1}^{n+1}-p_{1}^{*})=-(\frac{1}{2}(q_{1}^{n+1}+q_{1}^{*})+(q_{1}^{n+1}+q_{1}^{*})q_{2}^{n+1}), \\
&\frac{2}{\tau}(p_{2}^{n+1}-p_{2}^{*})=
-(\frac{1}{2}(q_{2}^{n+1}+q_{2}^{*})+(q_{1}^{*})^{2})+\frac{1}{3}((q_{2}^{n+1})^{2}+q_{2}^{n+1}q_{2}^{*}+(q_{2}^{*})^{2}),
\end{aligned}
\end{equation}
and the corresponding PAVF-P method
\begin{equation}\label{eq:3-5}
\begin{aligned}
&\frac{1}{\tau}(q_{1}^{n+1}-q_{1}^{n})=\frac{1}{2}(p_{1}^{n+1}+p_{1}^{n}), \\
&\frac{1}{\tau}(q_{2}^{n+1}-q_{2}^{n})=\frac{1}{2}(p_{2}^{n+1}+p_{2}^{n}), \\
&\frac{1}{\tau}(p_{1}^{n+1}-p_{1}^{n})=-\frac{1}{2}(q_{1}^{n+1}+q_{1}^{n}+(q_{1}^{n+1}+q_{1}^{n})(q_{2}^{n+1}+q_{2}^{n})), \\
&\frac{1}{\tau}(p_{2}^{n+1}-p_{2}^{n})=
-\frac{1}{2}(q_{2}^{n+1}+q_{2}^{n}+(q_{1}^{n+1})^{2}+(q_{1}^{n})^{2})+\frac{1}{3}((q_{2}^{n+1})^{2}+q_{2}^{n+1}q_{2}^{n}+(q_{2}^{n})^{2}),
\end{aligned}
\end{equation}
for the H\'{e}non-Heiles system \eqref{eq:3-1}. As aforementioned, both the two methods are energy-preserving and of second order accuracy. In addition, the PAVF-C method can inherit the semi-implicit property which however is no long satisfied by the PAVF-P method in view of the third equation in \eqref{eq:3-5}.

\subsubsection{Numerical experiments}

In the following experiments, we will test two classic orbits of the H\'{e}non-Heiles system \eqref{eq:3-1}: (i) chaotic orbits $(H^0=1/6)$, (ii)  box orbits $(H^0=0.02)$ with initial conditions $(q_1,q_2,p_2)=(0.1,-0.5,0)$ and $(0,-0.082,0)$. The leaving $p_1$ is found from the energy function \eqref{eq:3-1}. We solve the two cases by the AVF method \eqref{eq:3-2}, the PAVF method \eqref{eq:3-3},
the PAVF-C method \eqref{eq:3-4} and the PAVF-P method \eqref{eq:3-5}, respectively. The temporal increment is always set to $\tau=0.2$. The relative energy error is defined by
\begin{equation*}
RH^{n}=\big|(H^{n}-H^{0})/H^{0}\big|,
\end{equation*}
where $H^{n}$ denotes the Hamiltonian at $t=t^{n}$.\\

\begin{figure}[H]
	\centering
	(a)\includegraphics[width=0.45\textwidth]{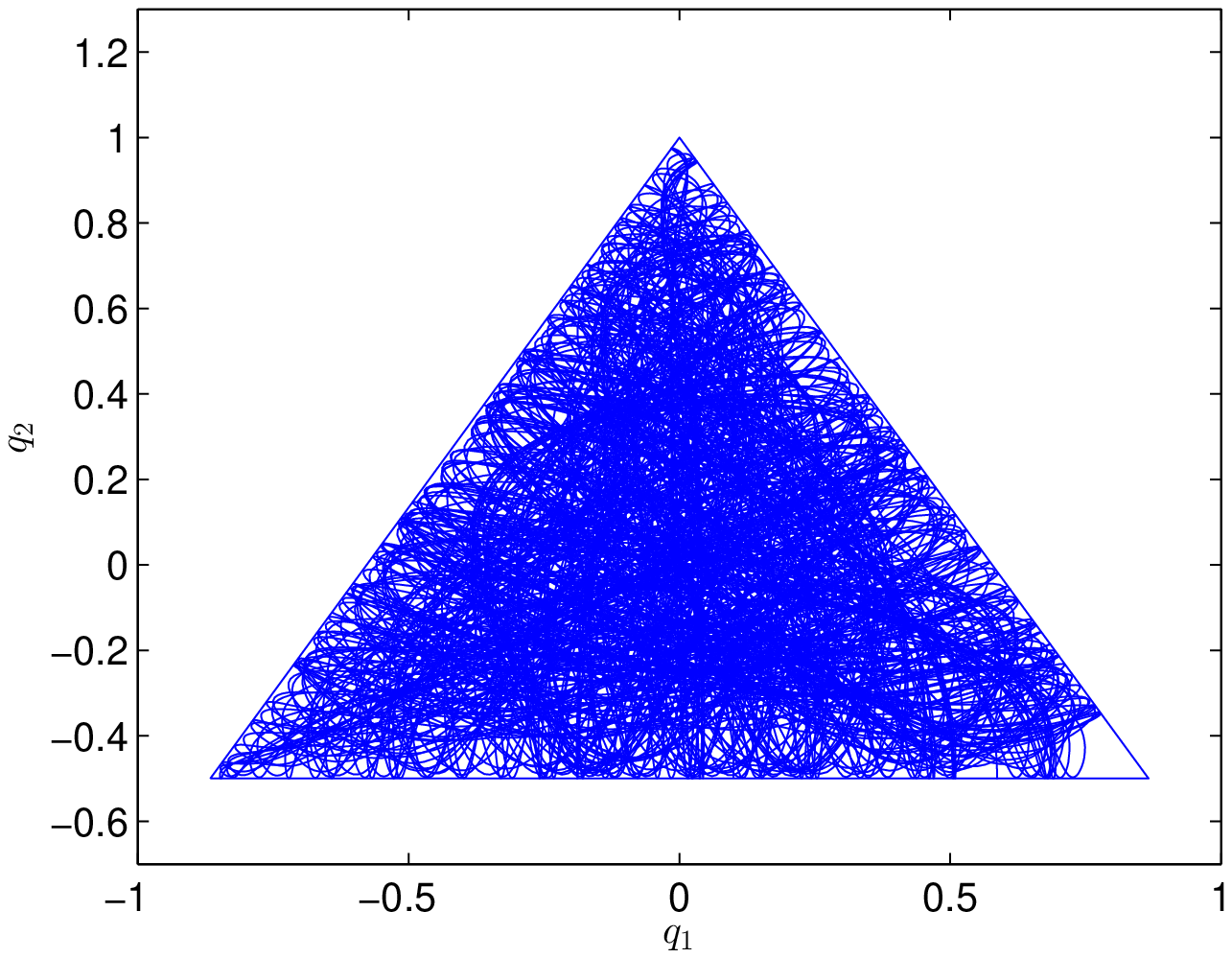}
	(b)\includegraphics[width=0.45\textwidth]{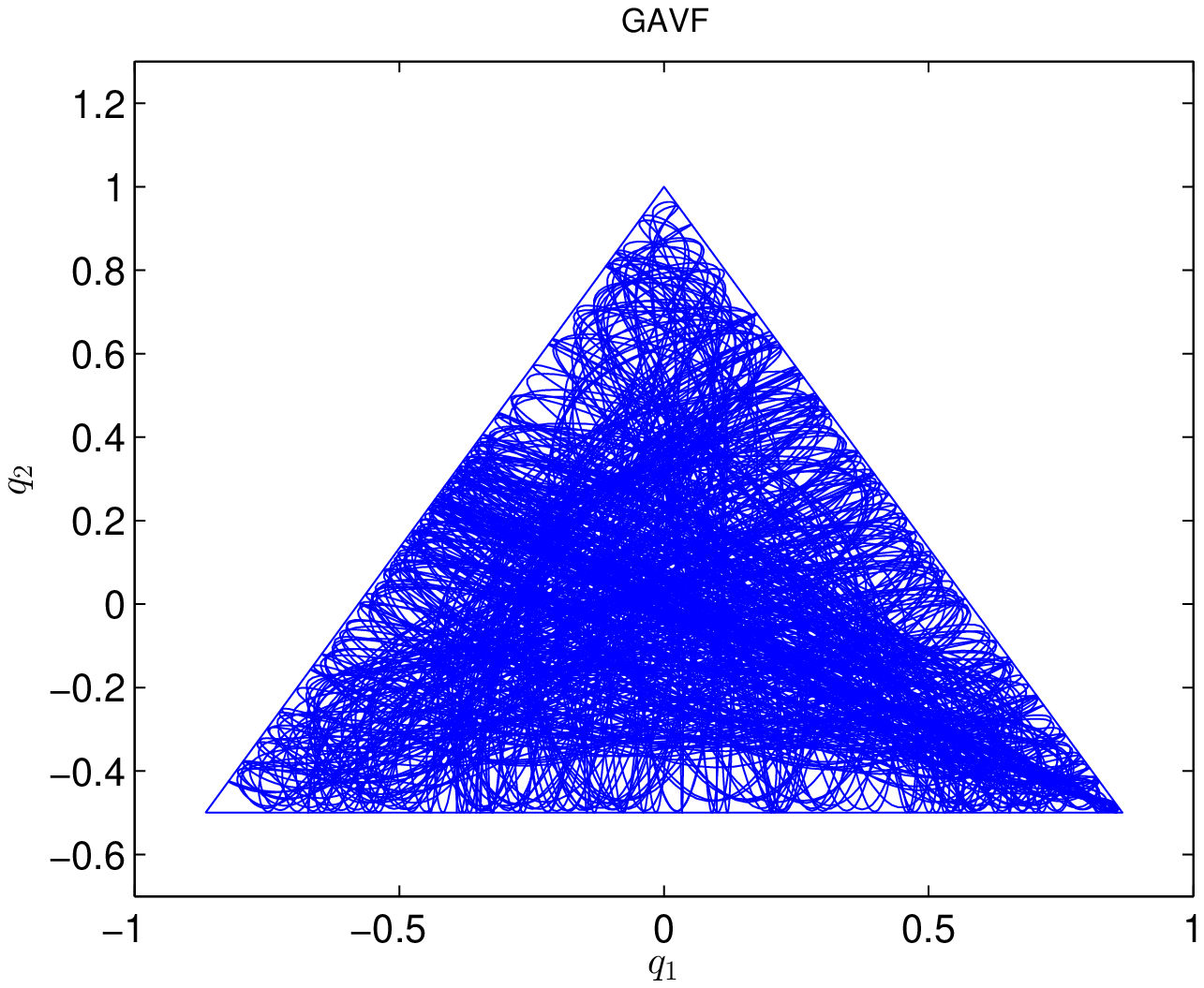}\\
	(c)\includegraphics[width=0.45\textwidth]{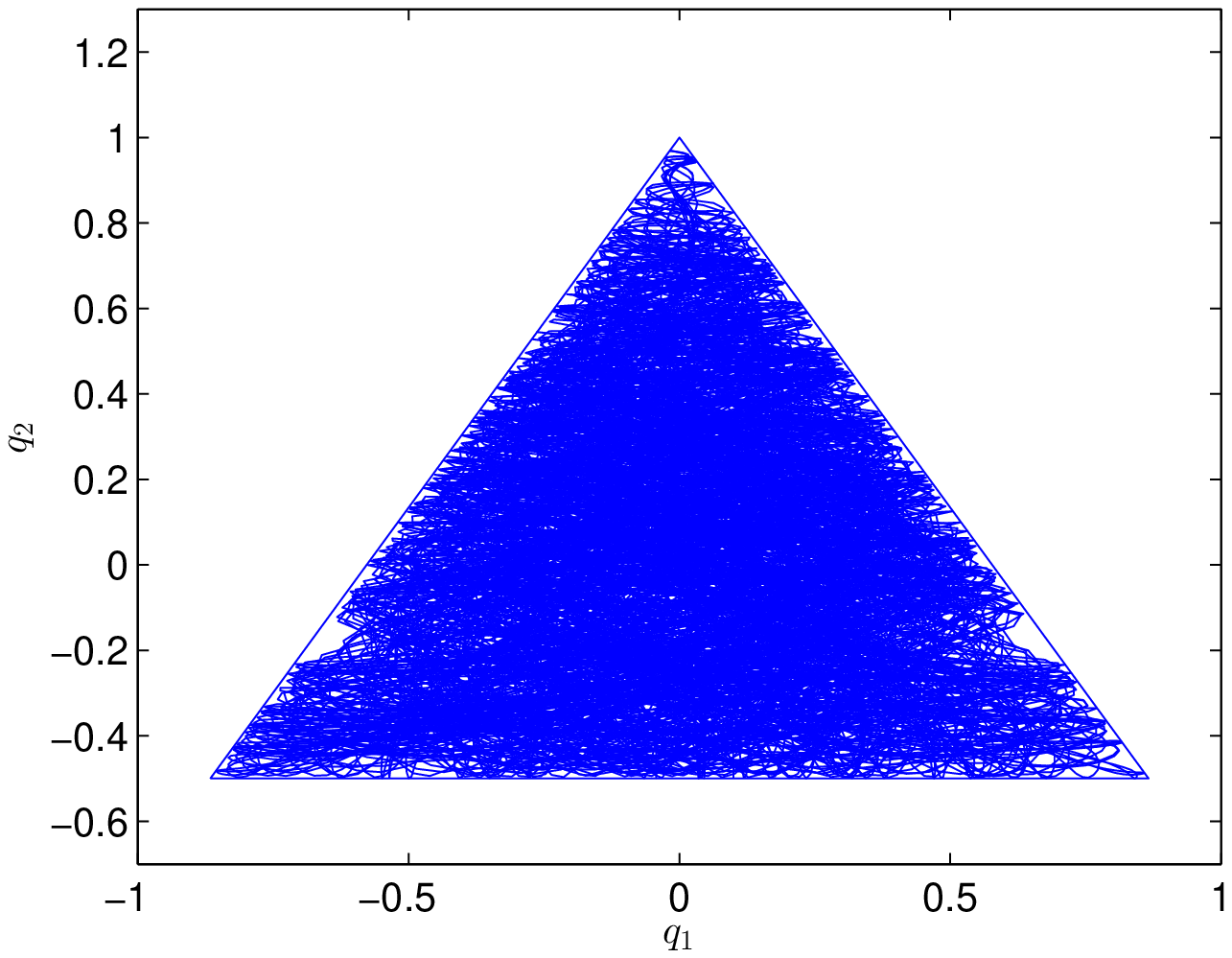}
	(d)\includegraphics[width=0.45\textwidth]{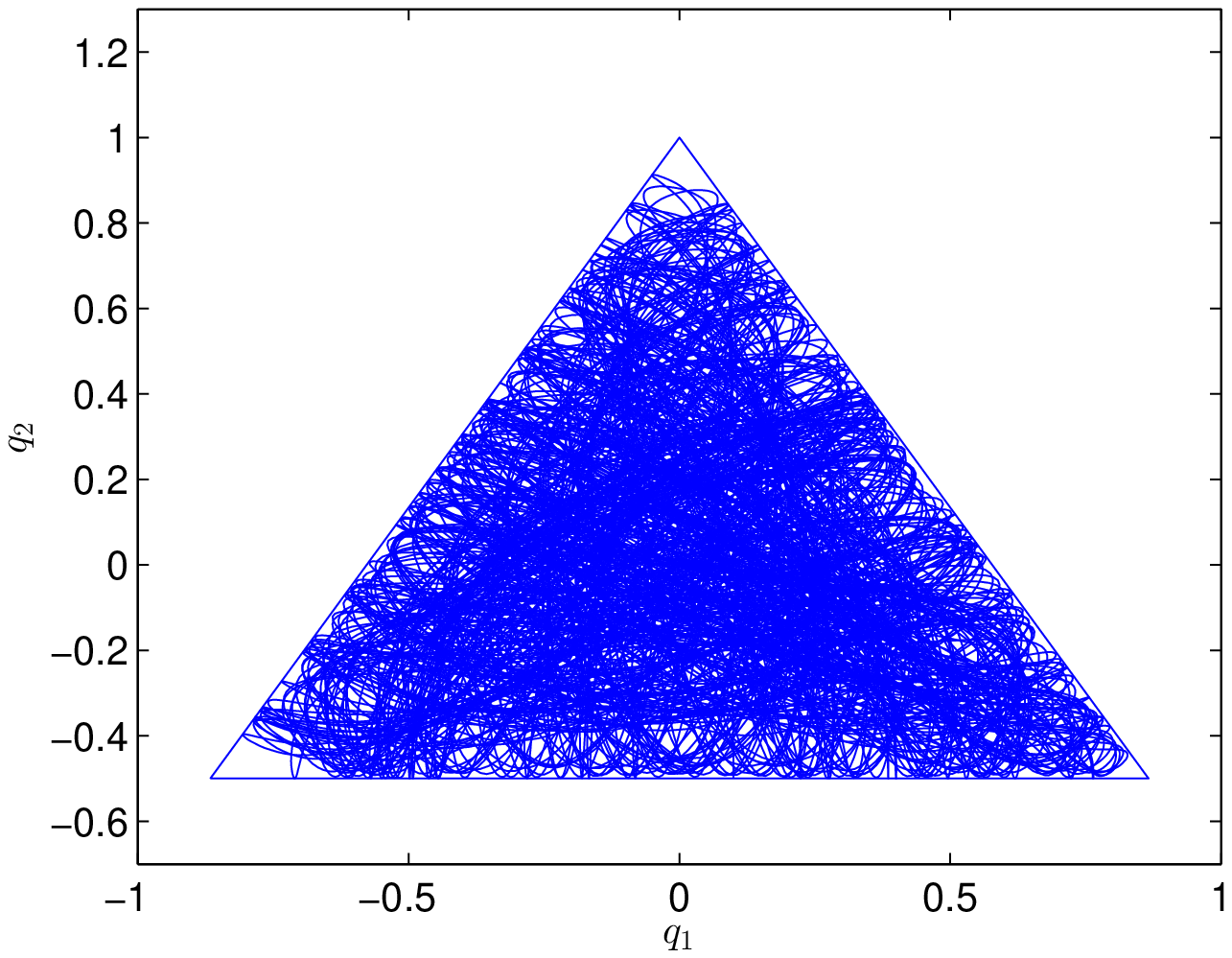}
	\caption{Chaotic orbits of the four energy-preserving methods for the H\'{e}non-Heiles system till $t=2000$. (a): The AVF method; (b): The PAVF method; (c): The PAVF-C method; (d): The PAVF-P method.}
	\label{fig:3-1}
\end{figure}

\begin{figure}[H]
	\centering
	(a)\includegraphics[width=0.45\textwidth]{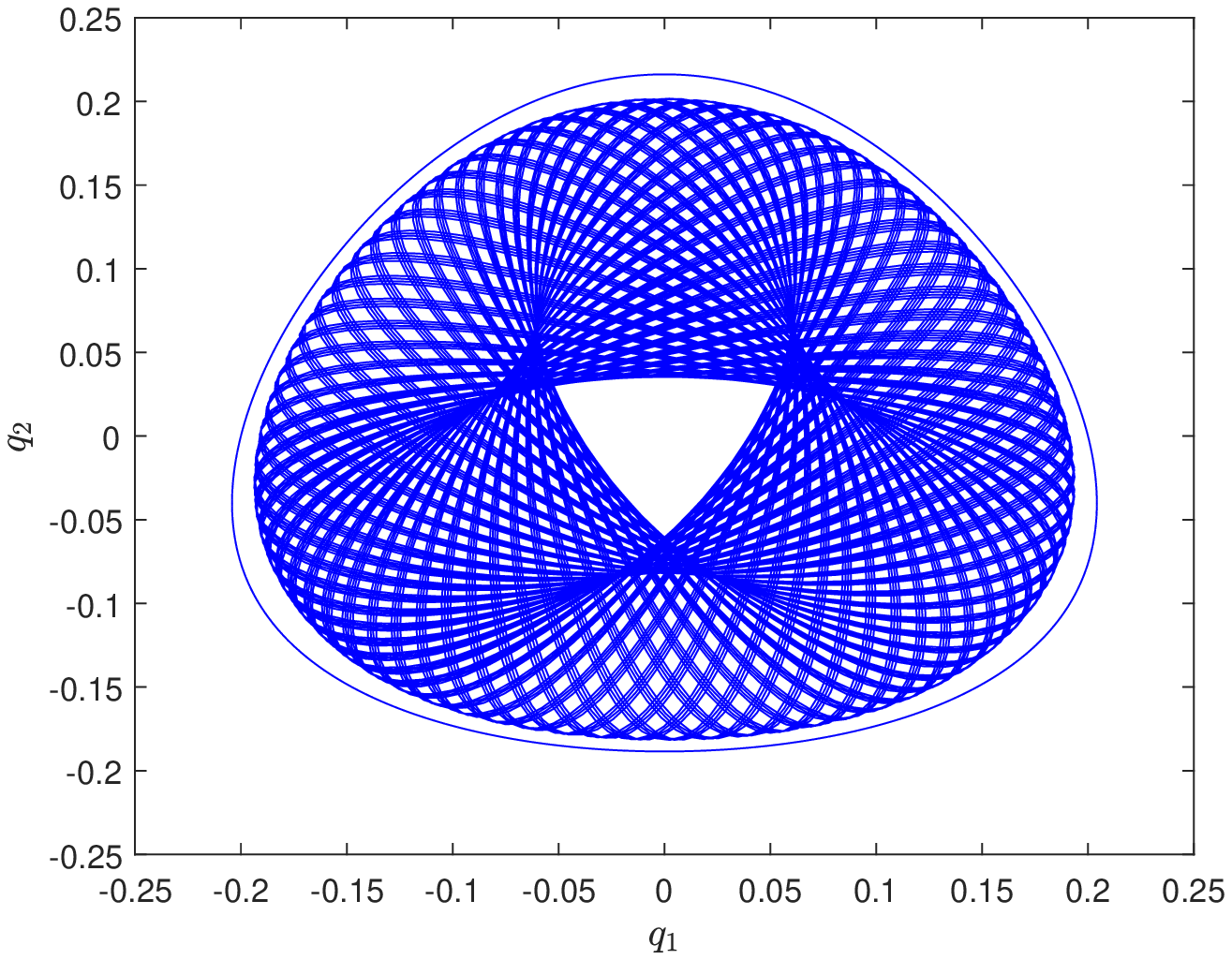}
	(b)\includegraphics[width=0.45\textwidth]{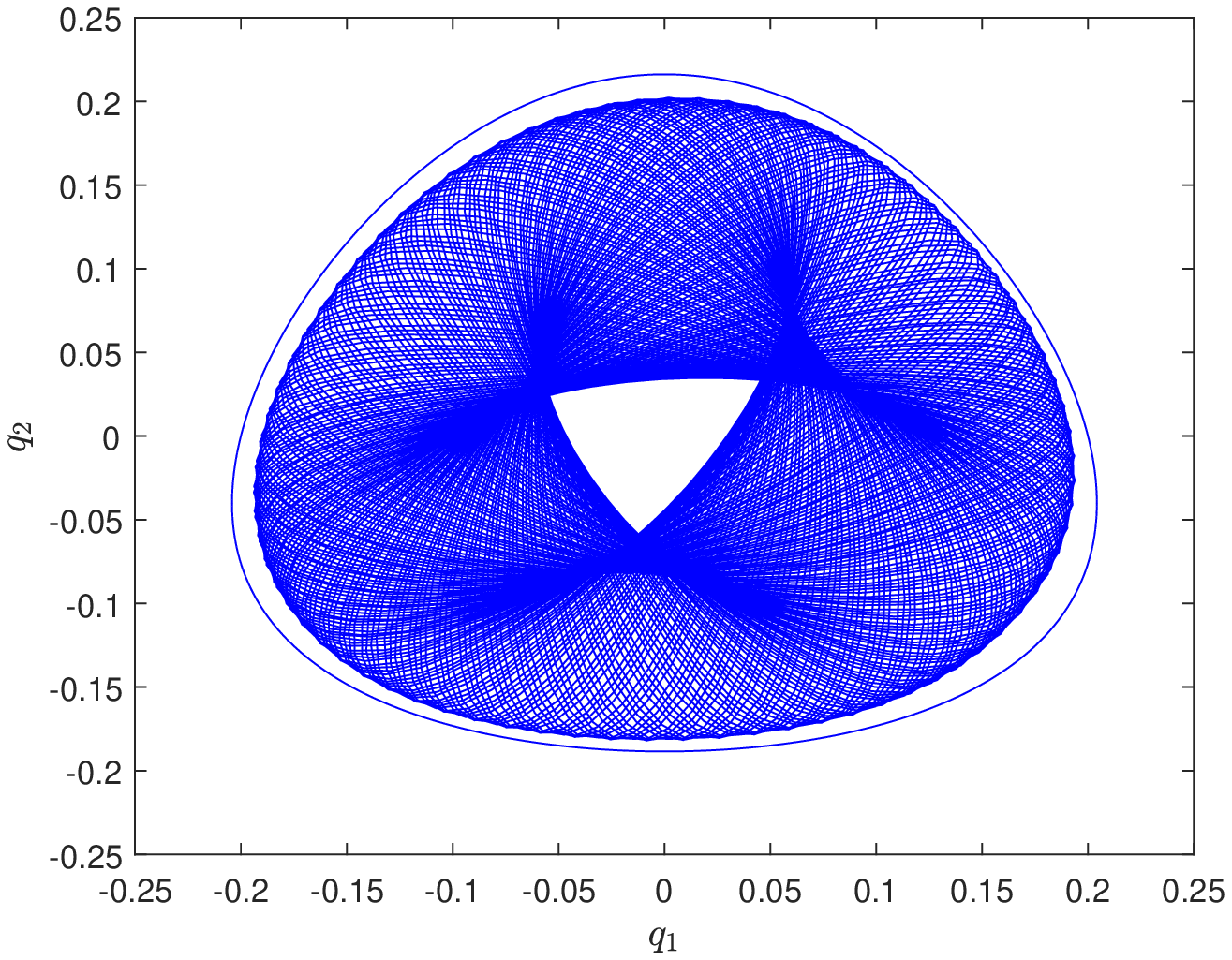}\\
	(c)\includegraphics[width=0.45\textwidth]{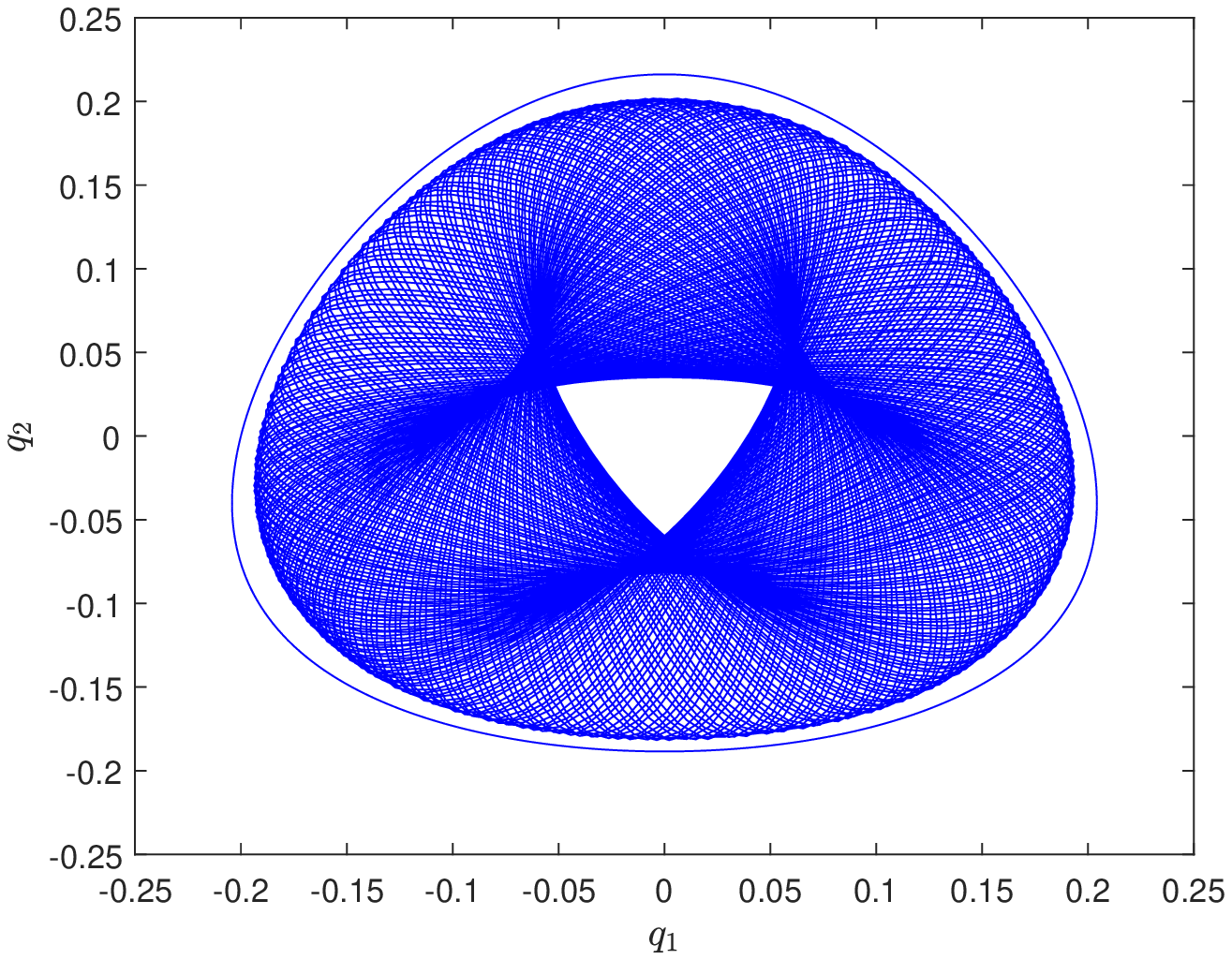}
	(d)\includegraphics[width=0.45\textwidth]{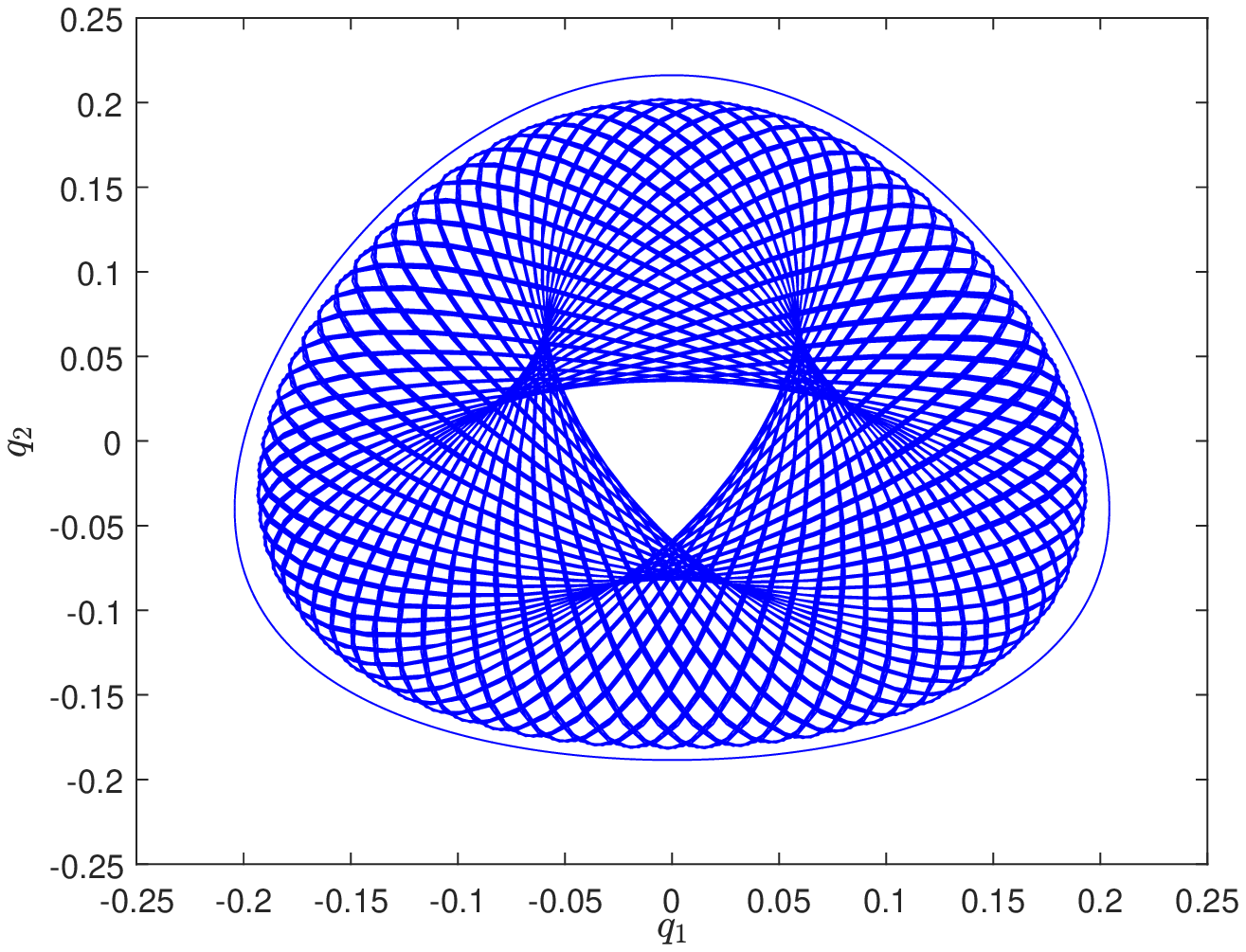}
	\caption{Box orbits of the four energy-preserving methods for the H\'{e}non-Heiles system till $t=1000$. (a): The AVF method; (b): The PAVF method; (c): The PAVF-C method; (d): The PAVF-P method.}
	\label{fig:3-11}
\end{figure}

\begin{figure}[H]
	\centering
	(a)\includegraphics[width=0.45\textwidth]{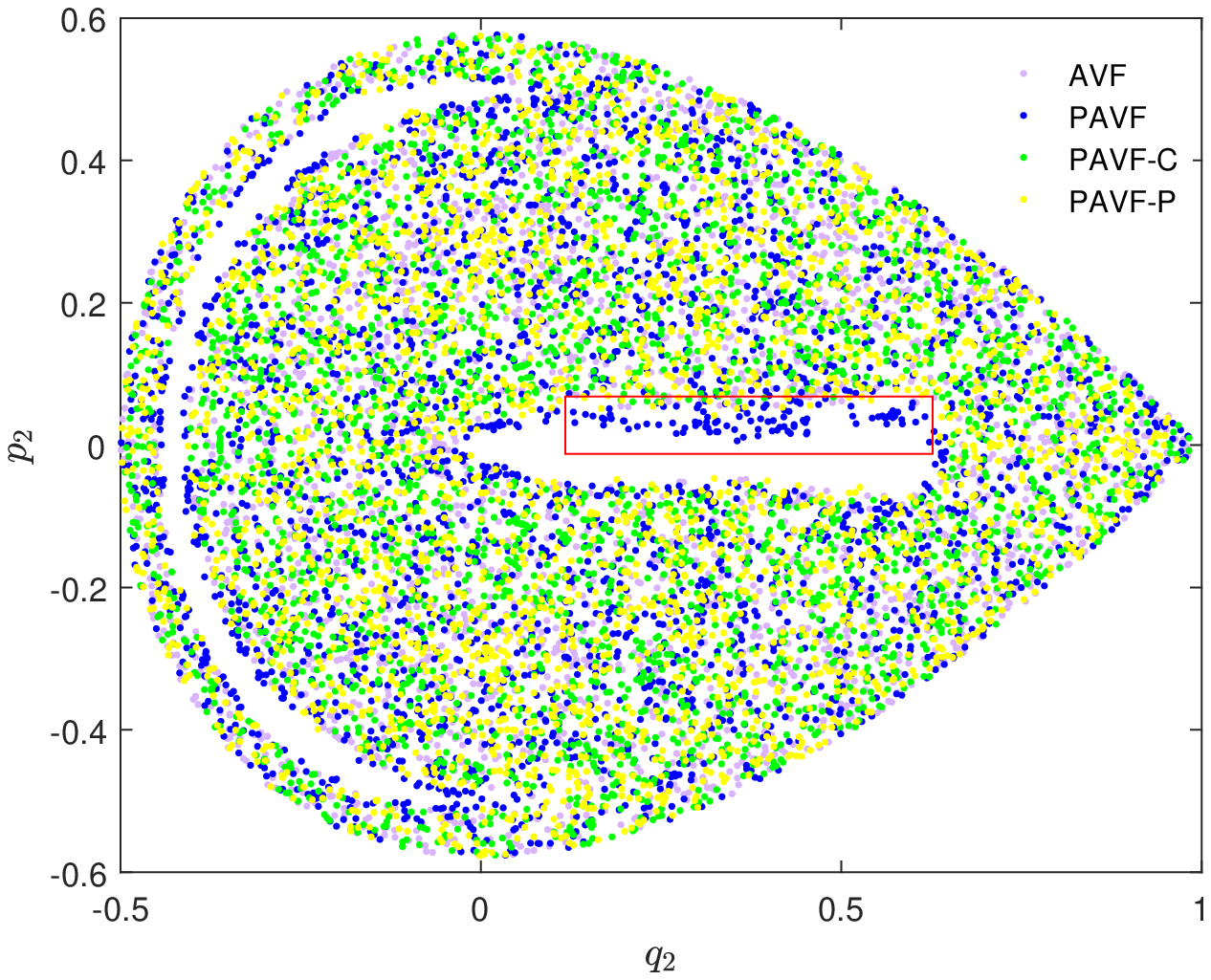}
	(b)\includegraphics[width=0.45\textwidth]{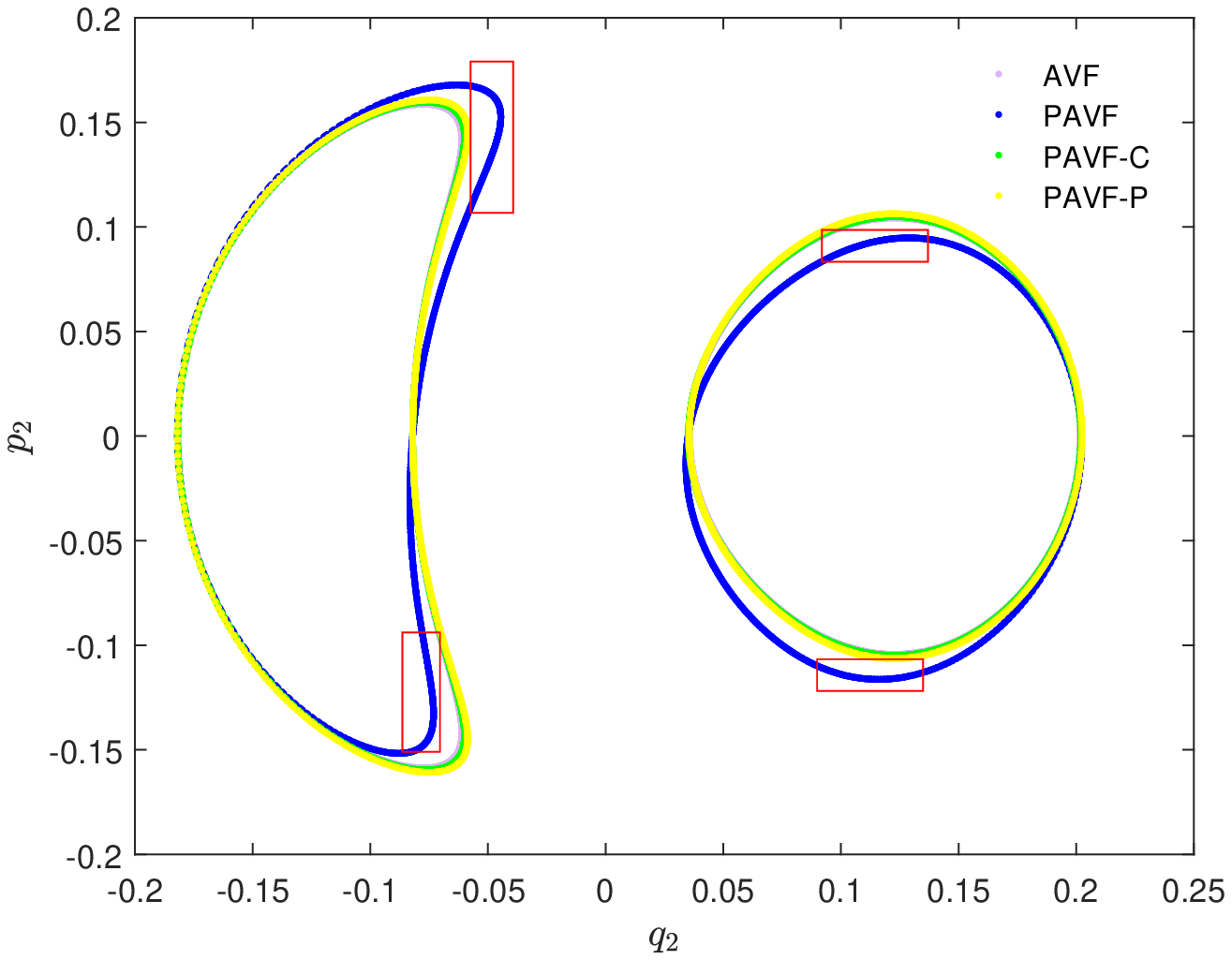}
	\caption{Poincar\'e cuts for the four numerical method (a): $H^0=1/6$; (b): $H^0=0.02$. The blue dots in red boxes show a little drift of obits produced by the first-order PAVF method.}
	\label{fig:3-11a}
\end{figure}

\begin{figure}[H]
	\centering
	(a)\includegraphics[width=0.45\textwidth]{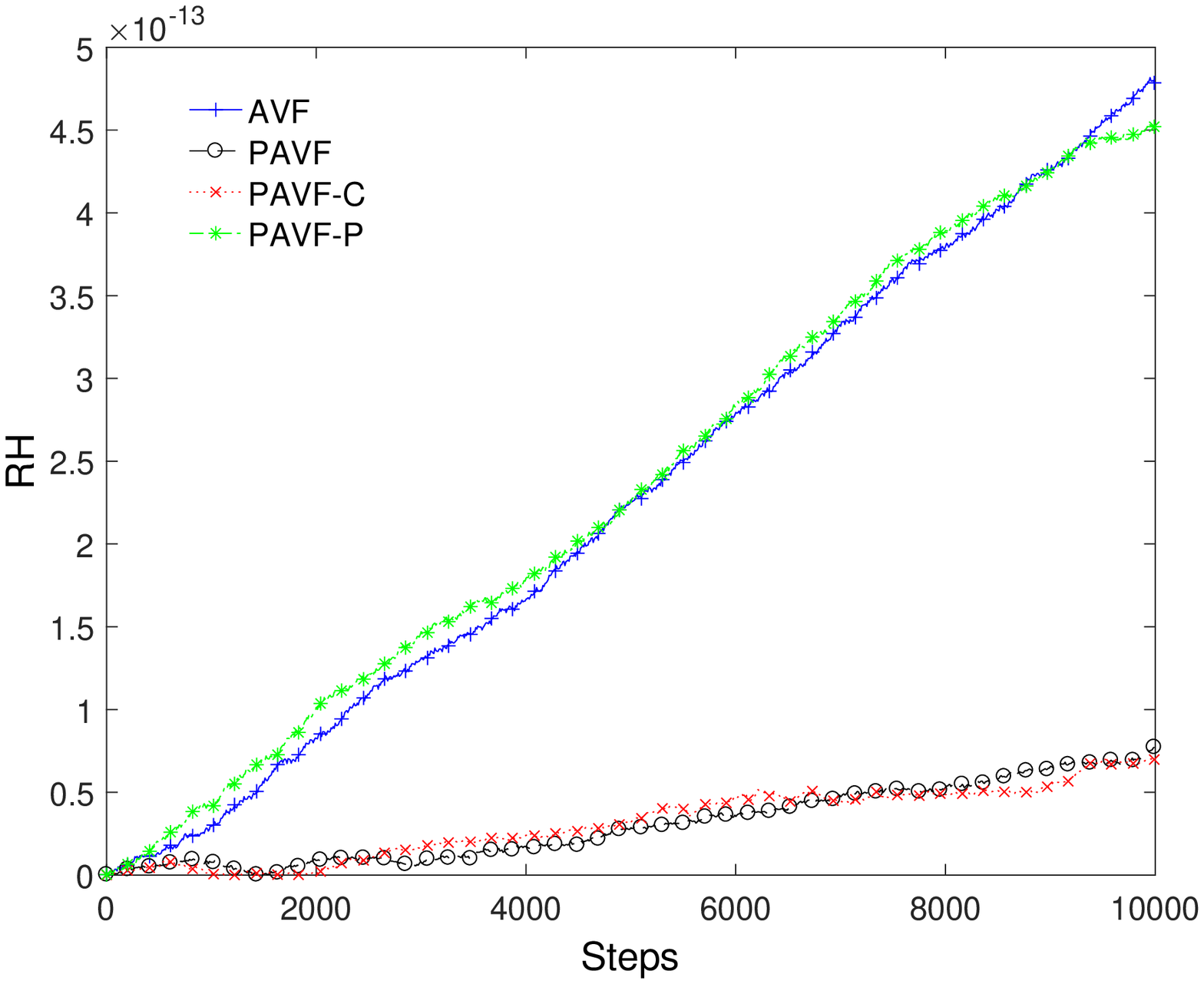}
	(b)	\includegraphics[width=0.45\textwidth]{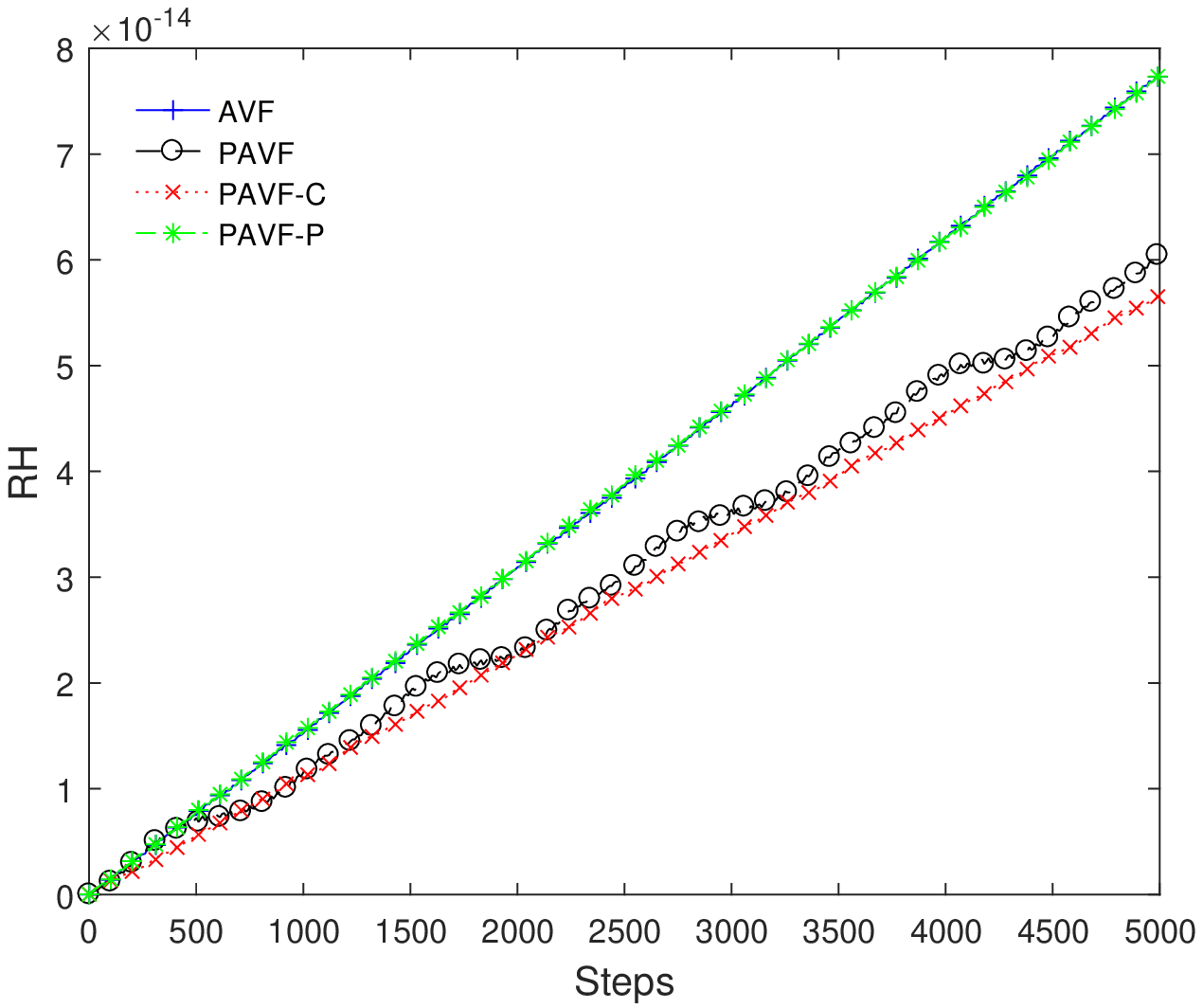}
	\caption{Relative energy error of the four energy-preserving methods for the H\'{e}non-Heiles system with chaotic orbits (a) and box orbits (b).}
	\label{fig:3-2}
\end{figure}

In Figure.~\ref{fig:3-1} it is clear that under the critical energy $1/6$ all the methods can ensure the numerical trajectory never escape the triangle for any value of $t$ and present the chaotic orbits as expected. While for the second case $H^0=0.02$, the four methods can also recurrence the box orbits in Figure.~\ref{fig:3-11}. Besides the result from the PAVF method, the inner shape of orbits all exhibits like an equilateral triangle. While the associate shape of the PAVF method suffers a little deformation which is probably caused by its lower accuracy. This can be also conformed from the corresponding Poincar\'e cuts in Figure.~\ref{fig:3-11a}. Furthermore, we find that the density between the orbits are very relied on the temporal increment. If we set $\tau=0.21$, the orbits of the PAVF-P method are very similar to that of the PAVF-C method. In Figure.~\ref{fig:3-2}, we present the corresponding errors in the discrete energy of all methods which uniformly reach the machine accuracy but  with a linear growth mainly caused by the numerical iteration. In order to compare the computational cost more clearly, we run the simulation for a longer time till $t=2\times 10^6$ for both cases. From Table.~\ref{tab:3-1}, we can see that the PAVF and PAVF-C schemes have significant advantages in CPU time than the rest ones. In addition, due to the fully implicit property of the AVF scheme and PAVF-P scheme, their CPU time is comparable to each other but more expensive than the semi-implicit PAVF and PAVF-C schemes.

\begin{table}[H]
	\centering
\caption{Computational cost for the H\'{e}non-Heiles system by four energy-preserving methods.}
\label{tab:3-1}
\begin{tabular*}{\textwidth}{@{\extracolsep{\fill}}ccccc}
	\hline
 $t=2\times 10^6$    & AVF                      & PAVF                      & PAVF-C                     & PAVF-P \\   \hline\noalign{\smallskip}
Case 1     &      16.8125            &     1.2182              & 2.4219                   & 16.6250\\\hline
Case 2     &      15.0781            &     1.1566              & 2.3436                   & 14.6563\\\hline
\end{tabular*}
\end{table}

\subsection{Klein-Gordon-Schr\"{o}dinger equation}
\label{sec:3-2}

We consider the following Klein-Gordon-Schr\"{o}dinger (KGS) equation:
\begin{equation}\label{eq:3-6}
\begin{aligned}
& i \varphi_{t}+ \frac{1}{2}\varphi_{xx}+ u \varphi = 0,\\ 
& u_{tt}- u_{xx}+u- |\varphi|^{2}= 0,\\ 
\end{aligned} \quad x\in (x_{L},x_{R}),~ t>0, ~i=\sqrt{-1},
\end{equation}
which describes a system of conserved scalar nucleons interacting with neutral scalar Mesons coupled with Yukawa interaction, where $\phi(x,t)$ represents a complex scalar nucleon field, $u(x,t)$ a real scalar meson field \cite{hjl09}.

The conservation of energy is a crucial property for the KGS equations.
Let $\varphi=q+pi$ and $v=\frac{1}{2}u_{t}.$ We can rewrite \eqref{eq:3-6} as a first-order system
\begin{equation}\label{eq:3-7}
\begin{aligned}
& u_{t} = 2v,\\ 
& v_{t} = \frac{1}{2}(u_{xx}-u+ (p^{2}+q^{2}) ),\\
& p_{t} =\frac{1}{2} q_{xx} + u q,\\
& q_{t} =- \frac{1}{2} p_{xx} - u p.
\end{aligned}
\end{equation}
Under homogeneous Dirichlet boundary conditions
\begin{align*}
& p(x_{L},t)=p(x_{R},t)=0,\quad q(x_{L},t)=q(x_{R},t)=0,\\
& u(x_{L},t)=u(x_{R},t)=0,\quad v(x_{L},t)=v(x_{R},t)=0,
\end{align*}
the above system \eqref{eq:3-7} can be comprised to an infinite-dimensional Hamiltonian
system
\begin{equation}\label{eq:3-8}
\frac{dz}{dt}=S \frac{\delta \mathcal{H}(z)}{\delta z},
\end{equation}
where $z=(u,v,p,q)^{T}$, 
$
S=
\left[
\begin{array}{rrrr}
0   & 1 & 0   & 0  \\
-1   & 0 & 0   & 0  \\
0   & 0 & 0   & -1  \\
0   & 0 & 1   & 0
\end{array}
\right],
$
and the Hamiltonian functional is defined by
\begin{equation*}
H(z)= \int_{x_{L}}^{x_{R}}\frac{1}{4}\left[p_{x}^{2}+q_{x}^{2}+u^{2}+u_{x}^{2}+4v^{2}-2u(p^{2}+q^{2})\right]\rd x.
\end{equation*}
Besides the energy conservation law $H(z(t))=H(z(0))$, the KGS system further obeys the mass conservation law
\begin{equation*}
M(z(t))=\int_{x_{L}}^{x_{R}}(p^{2}+q^{2})dx=M(z(0)).
\end{equation*}

\subsubsection{Derivation of the partitioned AVF schemes}

Notice that in the current work we mainly focus on the partitioned AVF methods for Hamiltonian systems, therefore, we just take the conventional spatial discretization for illustration. Hereafter, we choose the second-order central difference operator to approximate the spatial derivatives in \eqref{eq:3-7}. To present the concrete numerical schemes, we first introduce some notations.

The spatial-temporal domain is discretized as follows: $t_{n}=n \tau$, for $n=0,1,...$, and $x_{j}=x_{L}+jh$, for $j=0,1,...,J$ where $h=L/J$ and $\tau$ are the spatial and temporal steps. Let $(u_{j}^{n},v_{j}^{n},p_{j}^{n},q_{j}^{n})$
be the numerical approximations to the exact solutions $(u(x,t)$, $v(x,t)$, $p(x,t)$, $q(x,t))$ at the grid point $(x_{j},t_{n})$. The corresponding vector forms at any time level are then denoted by
\[
\begin{aligned}
&U=(u_{0},u_{1},\ldots,u_{J})^{T},~~V=(v_{0},v_{1},\ldots,v_{J})^{T},\\ &P=(p_{0},p_{1},\ldots,p_{J})^{T},~~Q=(q_{0},q_{1},\ldots,q_{J})^{T}.
\end{aligned}
\]
With these notations, we can define the inner product as well as norms of vectors:
\[
(U,V)=h\sum_{j=0}^{J}u_j\bar{v}_j, \quad \|U\|_J=(U,U)^{\frac{1}{2}},\quad \|U\|_\infty=\max_{j}|u_j|.
\]

The semi-discretization of the KGS equation \eqref{eq:3-7} by the central difference scheme can be written as
\begin{equation}\label{eq:3-9}
\begin{aligned}
&U_{t}=2V,\\ 
&V_{t}=\frac{1}{2}DU- \frac{1}{2}U+\frac{1}{2}(P^{2}+Q^{2}),\\ 
&P_{t}=\frac{1}{2}DQ+ U\cdot Q,\\
&Q_{t}=-\frac{1}{2}DP- U\cdot P,
\end{aligned}
\end{equation}
where $D$ represents the central differentiation matrix, and the operation $P^{2}=P\cdot P$.
Here, $`\cdot$' means the point multiplication between vectors, that is, $P\cdot Q=(p_{0}q_{0},p_{1}q_{1},\ldots,p_{J}q_{J})^{T}$. In analogy to the continuous case, we can also transform the system \eqref{eq:3-9} back to a finite-dimensional canonical Hamiltonian form
\begin{equation}\label{eq:3-10}
\frac{dZ}{dt}=f(Z)=S \nabla H(Z),
\end{equation}
where $Z=(U^{T},V^{T},P^{T},Q^{T})^{T}$, $
S=
\left[
\begin{array}{rrrr}
0   & I & 0   & 0  \\
-I   & 0 & 0   & 0  \\
0   & 0 & 0   & -I  \\
0   & 0 & I   & 0
\end{array}
\right],
$
and the discrete Hamiltonian is defined by
\[
H(Z)=\frac{1}{4}\left(-P^{T}
DP-Q^{T}DQ-U^{T}DU
+U^{T}U+4V^{T}V
-2\sum_{j=0}^{J}u_{j}\left(p_{j}^{2}+q_{j}^{2}\right)\right).
\]
From system \eqref{eq:3-10}, we can derive the semi-discrete energy conservation law
\[
\frac{dH(Z(t))}{dt}=\nabla H(Z)^{T}f(Z)=\nabla H(Z)^{T}S\nabla H(Z)=0,
\]
and letting $M(Z)=\| P\|_{J}^{2}+\| Q\|_{J}^{2}$, we also have the semi-discrete mass conservation law
\[
\frac{dM(Z)}{dt}  =  2h\left(\frac{1}{2}P^{T}
DQ-\frac{1}{2}Q^{T}DP+P^{T}
(U\cdot Q)-Q^{T}(U\cdot P)\right)=0.
\]
Hence, the flow of the semi-discrete system \eqref{eq:3-9} preserves both the total energy
$H(Z)$ and mass $M(Z)$ exactly.

The following part is devoted to the construction of the partitioned AVF methods \eqref{eq:2-7}, \eqref{eq:2-9} and \eqref{eq:2-10} for the KGS equation \eqref{eq:3-6}. For comparison we also present the original second-order AVF scheme 
\begin{equation}\label{eq:3-11}
\begin{aligned}
&\delta_{t}^{+}U^{n}= V^{n}+V^{n+1},\\ 
&\delta_{t}^{+}V^{n}=\frac{1}{2}(D U^{n+\frac{1}{2}}- U^{n+\frac{1}{2}})+\frac{1}{12}((P^{n+1})^{2}+4(P^{n+\frac{1}{2}})^{2}+(P^{n})^{2}\\
&\hspace{0.4\linewidth}+(Q^{n+1})^{2}+4(Q^{n+\frac{1}{2}})^{2}+(Q^{n})^{2}),\\ 
&\delta_{t}^{+}P^{n}=\frac{1}{2}D Q^{n+\frac{1}{2}}+ \frac{1}{6}(U^{n+1}\cdot Q^{n+1}+4U^{n+\frac{1}{2}}\cdot Q^{n+\frac{1}{2}}
+U^{n}\cdot Q^{n}),\\ 
&\delta_{t}^{+}Q^{n}=-\frac{1}{2}D P^{n+\frac{1}{2}}- \frac{1}{6}(U^{n+1}\cdot P^{n+1}+4U^{n+\frac{1}{2}}\cdot P^{n+\frac{1}{2}}
+U^{n}\cdot P^{n}),
\end{aligned}
\end{equation}
where $\delta_t^+$ is the standard forward differential operator. Obviously, the AVF method for the KGS equation is a fully implicit scheme but with the energy-preserving property.

Next, we introduce our partitioned AVF methods. Under this circumstance, the dimension parameters in  \eqref{eq:2-6} are set to $\tilde{m}=4$ and $d=J+1$. Then the corresponding PAVF scheme yields
\begin{equation}\label{eq:3-12}
\begin{aligned}
&\delta_{t}^{+}U^{n}=\int_{0}^{1}H_{V}(U^{n+1},\xi V^{n+1}+(1-\xi)V^{n},P^{n},Q^{n})d\xi,\\
&\delta_{t}^{+}V^{n}=-\int_{0}^{1}H_{U}(\xi U^{n+1}+(1-\xi)U^{n},V^{n},P^{n},Q^{n})d\xi,\\ 
&\delta_{t}^{+}Q^{n}=-\int_{0}^{1}H_{Q}(U^{n+1},V^{n+1},P^{n+1},\xi Q^{n+1}+(1-\xi)Q^{n})d\xi,\\ 
&\delta_{t}^{+}Q^{n}=\int_{0}^{1}H_{P}(U^{n+1},V^{n+1},\xi P^{n+1}+(1-\xi)P^{n},Q^{n})d\xi.
\end{aligned}
\end{equation}
which can be further integrated as
\begin{equation}\label{eq:3-13}
\begin{aligned}
&\delta_{t}^{+}U^{n}=V^{n}+V^{n+1},\\
&\delta_{t}^{+}V^{n}=\frac{1}{2}(D U^{n+\frac{1}{2}}- U^{n+\frac{1}{2}}+(P^{n})^{2}+(Q^{n})^{2}),\\ 
&\delta_{t}^{+}P^{n}=\frac{1}{2}D Q^{n+\frac{1}{2}}+ U^{n+1}\cdot Q^{n+\frac{1}{2}},\\ 
&\delta_{t}^{+}Q^{n}=-\frac{1}{2}D P^{n+\frac{1}{2}}- U^{n+1}\cdot P^{n+\frac{1}{2}}.
\end{aligned}
\end{equation}
Apparently, the PAVF scheme \eqref{eq:3-13} is simpler than the AVF scheme \eqref{eq:3-11} that just requires to solve two sets of linear algebraic equations in contrast to apply numerical iterations for the entire system.

\begin{theorem}\label{thm:3-1}
	The PAVF method \eqref{eq:3-13} is energy and mass conservative. That is,
	its solution satisfies the following conservation laws
	\begin{equation*}
	\begin{aligned}
	&M(U^{n},V^{n},P^{n},Q^{n})= M(U^{0},V^{0},P^{0},Q^{0}),\\
	 &H(U^{n},V^{n},P^{n},Q^{n})= H(U^{0},V^{0},P^{0},Q^{0}),
	\end{aligned}\quad \forall ~n=1,2,....
	\end{equation*}
\end{theorem}
\begin{proof}
	Taking the discrete inner products of the last two equations in \eqref{eq:3-13} with $P^{n+1}+P^{n}$ and $Q^{n+1}+Q^{n}$ respectively, we obtain
	\begin{equation}
	\label{eq:3-14}
	\frac{h}{\tau}(P^{n+1}+P^{n})^{T}(P^{n+1}-P^{n})=h(P^{n+1}+P^{n})^{T}(\frac{1}{2}D Q^{n+\frac{1}{2}}+ U^{n+1}\cdot Q^{n+\frac{1}{2}}),
	\end{equation}
	and
	\begin{equation}
	\label{eq:3-15}
	\frac{h}{\tau}(Q^{n+1}+Q^{n})^{T}(Q^{n+1}-Q^{n})=h(Q^{n+1}+Q^{n})^{T}(-\frac{1}{2}D P^{n+\frac{1}{2}}- U^{n+1}\cdot P^{n+\frac{1}{2}}).
	\end{equation}
	Summing \eqref{eq:3-14} and \eqref{eq:3-15} together, we get the mass conservation law
	\begin{equation}
	\label{eq:5a2b30}
	\frac{1}{\tau}\left(\| P^{n+1}\|_{J}^{2}+\| Q^{n+1}\|_{J}^{2}-
	\| P^{n}\|_{J}^{2}-\| Q^{n}\|_{J}^{2}\right)=0.
	\end{equation}

	Taking the discrete inner products of \eqref{eq:3-13} with $\delta_{t}^{+}V^{n}$, $-\delta_{t}^{+}U^{n}$, $-\delta_{t}^{+}Q^{n}$ and $\delta_{t}^{+}P^{n}$ respectively, we have the following equalities:
	\begin{equation}\label{eq:3-18}
	\begin{aligned}
	&h(\delta_{t}^{+}V^{n})^{T}\delta_{t}^{+}U^{n}=\frac{h}{\tau}(V^{n+1}-V^{n})^{T} (V^{n+1}+V^{n}),\\
	-&h(\delta_{t}^{+}U^{n})^{T}\delta_{t}^{+}V^{n}=-\frac{h}{2\tau}(U^{n+1}-U^{n})^{T} (D U^{n+\frac{1}{2}}- U^{n+\frac{1}{2}}+(P^{n})^{2}+(Q^{n})^{2}),\\
	-&h(\delta_{t}^{+}Q^{n})^{T}\delta_{t}^{+}P^{n}=-\frac{h}{\tau}(Q^{n+1}-Q^{n})^{T} (\frac{1}{2}D Q^{n+\frac{1}{2}}+ U^{n+1}\cdot Q^{n+\frac{1}{2}}),\\
	&h(\delta_{t}^{+}P^{n})^{T}\delta_{t}^{+}Q^{n}=\frac{h}{\tau}(P^{n+1}-P^{n})^{T} (-\frac{1}{2}D P^{n+\frac{1}{2}}- U^{n+1}\cdot P^{n+\frac{1}{2}}).
	\end{aligned}
	\end{equation}
	Then the energy conservation law comes from summing the above four equations together and rearranging the terms with respect to the time levels $n$ and $n+1$.
\end{proof}

To derive the second-order partitioned AVF methods, including the PAVF-C and PAVF-P methods, we first present the adjoint of the PAVF scheme \eqref{eq:3-13} for the KGS equation:
\begin{equation}\label{eq:3-19}
\begin{aligned}
&\delta_{t}^{+}U^{n}= V^{n}+V^{n+1},\\
&\delta_{t}^{+}V^{n}=\frac{1}{2}(D U^{n+\frac{1}{2}}- U^{n+\frac{1}{2}}+(P^{n+1})^{2}+(Q^{n+1})^{2}),\\
&\delta_{t}^{+}P^{n}=\frac{1}{2}D Q^{n+\frac{1}{2}}+ U^{n}\cdot Q^{n+\frac{1}{2}},\\
&\delta_{t}^{+}Q^{n}=-\frac{1}{2}D P^{n+\frac{1}{2}}- U^{n}\cdot P^{n+\frac{1}{2}},
\end{aligned}
\end{equation}
which can be similarly proved as an energy-mass-preserving method either. With the adjoint scheme \eqref{eq:3-19}, we can directly write down the corresponding PAVF-C scheme
\begin{equation}\label{eq:3-20}
\begin{aligned}
&\frac{1}{\tau}(U^{*}-U^{n})= V^{n}+V^{*},\\
&\frac{1}{\tau}(V^{*}-V^{n})=\frac{1}{4}(D (U^{*}+U^{n})-  (U^{*}+U^{n})+2(P^{n})^{2}+2(Q^{n})^{2}),\\
&\frac{1}{\tau}(P^{*}-P^{n})=\frac{1}{4}D (Q^{*}+Q^{n})+ \frac{1}{2}U^{*}\cdot (Q^{*}+Q^{n}),\\
&\frac{1}{\tau}(Q^{*}-Q^{n})=-\frac{1}{4}D (P^{*}+P^{n})- \frac{1}{2}U^{*}\cdot (P^{*}+P^{n}),\\
&\frac{1}{\tau}(U^{n+1}-U^{*})= V^{n+1}+V^{*},\\
&\frac{1}{\tau}(V^{n+1}-V^{*})=\frac{1}{4}(D (U^{*}+U^{n+1})-  (U^{*}+U^{n+1})+2(P^{n+1})^{2}+2(Q^{n+1})^{2}),\\
&\frac{1}{\tau}(P^{n+1}-P^{*})=\frac{1}{4}D (Q^{*}+Q^{n+1})+ \frac{1}{2}U^{*}\cdot (Q^{*}+Q^{n+1}),\\
&\frac{1}{\tau}(Q^{n+1}-Q^{*})=-\frac{1}{4}D (P^{*}+P^{n+1})- \frac{1}{2}U^{*}\cdot (P^{*}+P^{n+1}),
\end{aligned}
\end{equation}
and the PAVF-P scheme 
\begin{equation}\label{eq:3-21}
\begin{aligned}
&\delta_{t}^{+}U^{n}= V^{n}+V^{n+1},\\ 
&\delta_{t}^{+}V^{n}=\frac{1}{2}(D U^{n+\frac{1}{2}}- U^{n+\frac{1}{2}})+\frac{1}{4}((P^{n})^{2}+(Q^{n})^{2}+(P^{n+1})^{2}+(Q^{n+1})^{2}),\\ 
&\delta_{t}^{+}P^{n}=\frac{1}{2}D Q^{n+\frac{1}{2}}+ U^{n+\frac{1}{2}}\cdot Q^{n+\frac{1}{2}},\\ 
&\delta_{t}^{+}Q^{n}=-\frac{1}{2}D P^{n+\frac{1}{2}}- U^{n+\frac{1}{2}}\cdot P^{n+\frac{1}{2}},
\end{aligned}
\end{equation}
which both are symmetric methods of second-order accuracy and can also preserve the discrete energy and mass conservation laws.

\begin{remark}
	In \cite{zhang16} authors also present two energy-preserving schemes for the KGS equation based on the coordinate increment discrete gradient method, which are equivalent to the adjoint and plus schemes here respectively. However, through the partitioned AVF strategy we can easily construct other energy-preserving schemes. Specifically, the resulting partitioned AVF composition scheme can not only conserve the discrete energy but also save computational time with second-order accuracy.
\end{remark}

\subsubsection{Numerical experiments}

In this section, we consider the numerical results of the four energy-preserving methods
for the KGS system \eqref{eq:3-6}. To measure the conservative properties, we define the
relative energy and mass errors at $t=t_{n}$ as
\[
RH^{n}=\big|(H^{n}-H^{0})/H^{0}\big|,\quad RM^{n}=\big|(M^{n}-M^{0})/M^{0}\big|.
\]
The discrete $L^{2}$-error and  $L^{\infty}$-error for solutions of the KGS equation are calculated by
\begin{align*}
& \mbox{\rm$L^{2}$-error}_{n}(h,\tau)=\| u(t_n)-U^{n}\|_{J}+\| v(t_n)-V^{n}\|_{J}+\| p(t_n)-P^{n}\|_{J}+\| q(t_n)-Q^{n}\|_{J},\\ 
& \mbox{\rm$L^{\infty}$-error}_{n}(h,\tau)=\| u(t_n)-U^{n}\|_{\infty}+\| v(t_n)-V^{n}\|_{\infty}+\| p(t_n)-P^{n}\|_{\infty}+\| q(t_n)-Q^{n}\|_{\infty}.
\end{align*}

The first example is about the propagation of one soliton, in which the initial conditions are taken as
\begin{align*}
&\varphi(x,0)=\frac{3\sqrt{2}}{4\sqrt{1-c^{2}}}\mathrm{sech}^{2}\Big(\frac{1}{2\sqrt{1-c^{2}}}(x-x_{0})\Big)\exp(icx),\\
&u(x,0)=\frac{3}{4(1-c^{2})}\mathrm{sech}^{2}\Big(\frac{1}{2\sqrt{1-c^{2}}}(x-x_{0})\Big),\\
&u_{t}(x,0)=\frac{3v}{4(1-c^{2})^{3/2}}\mathrm{sech}^{2}\Big(\frac{1}{2\sqrt{1-c^{2}}}(x-x_{0})\Big)\tanh\Big(\frac{1}{2\sqrt{1-c^{2}}}(x-x_{0})\Big),
\end{align*}
where $|c|<1$ is the propagating velocity of the wave and $x_{0}$ is the
initial phase.
The exact solution is given by
\begin{align*}
&\varphi(x,t)=\frac{3\sqrt{2}}{4\sqrt{1-c^{2}}}\mathrm{sech}^{2}\Big(\frac{1}{2\sqrt{1-c^{2}}}(x-ct-x_{0})\Big)\exp\Big(i\big(cx+\frac{1-c^{2}+c^{4}}{2(1-c^{2})}t\big)\Big),\\
&u(x,t)=\frac{3}{4(1-c^{2})}\mathrm{sech}^{2}\Big(\frac{1}{2\sqrt{1-c^{2}}}(x-ct-x_{0})\Big).
\end{align*}
We set the space interval $x\in [-10,10]$ with parameters $c=-0.8$, $x_{0}=0$ and the homogeneous Dirichlet boundary conditions. Table.~\ref{tab:3-2} and \ref{tab:3-3} give the temporal and spatial accuracy of the four energy-preserving methods. We can find that all the methods are both of second order in space and time except the PAVF method which is only first-order accuracy in time. The accuracy test validates the correctness of our methods.

\begin{table}[H]
	\centering
	\caption{Temporal accuracy of the four energy-preserving methods with $h=0.02$.}\label{tab:3-2}
	\begin{tabular*}{\textwidth}[h]{@{\extracolsep{\fill}}l l l l l l} \hline
		& $\tau$ & $L^\infty$-error & order & $L^2$-error & order  \\ \hline
		\multirow{4}{*}{AVF} & 1/10  & 1.15E-02 & - & 1.28E-02 & - \\
		& 1/11 & 9.44E-03 & 2.04 & 1.05E-02 & 2.03  \\
		& 1/12 & 7.90E-03 & 2.05 & 8.81E-03 & 2.04  \\
		& 1/13 & 6.70E-03 & 2.06 & 7.48E-03 & 2.05 \\ [1ex] \hline		
		\multirow{4}{*}{PAVF} & 1/10  & 1.06E-01 & - & 1.44E-02 &  - \\
		& 1/11  & 9.64E-02 & 0.98 & 1.27E-01 & 1.00 \\
		& 1/12 & 8.85E-02 & 0.98 & 1.15E-01 & 1.00  \\
		& 1/13 & 8.18E-02 & 0.98 & 1.06E-01 & 1.00\\ [1ex]\hline			
		\multirow{4}{*}{PAVF-C} & 1/10  & 4.42E-03 & - & 4.19E-03 & - \\
		& 1/11  & 3.60E-03 & 2.14 & 3.42E-03 & 2.13 \\
		& 1/12 & 2.98E-03 & 2.17 & 2.84E-03 & 2.15 \\
		& 1/13 & 2.49E-03 & 2.21 & 2.38E-03 & 2.17 \\ [1ex]\hline		
		\multirow{4}{*}{PAVF-P} & 1/10  & 9.60E-03 & - & 1.07E-02 &  - \\
		& 1/11  & 7.89E-03 & 2.06 & 8.82E-03 & 2.02 \\
		& 1/12 & 6.59E-03 & 2.07 & 7.40E-03 & 2.03  \\
		& 1/13 & 5.57E-03 & 2.09 & 6.29E-03 & 2.03\\
		\hline
	\end{tabular*}
\end{table}

\begin{table}[H]
	\centering
	\caption{Spatial accuracy of the four energy-preserving methods with $\tau=0.001$.}\label{tab:3-3}
	\begin{tabular*}{\textwidth}[h]{@{\extracolsep{\fill}}c l l l l l} \hline
		& $h$ & $L^\infty$-error & order & $L^2$-error & order  \\ \hline
		\multirow{4}{*}{AVF}  & 2/10  & 6.04E-02 & - & 6.50E-02 & - \\
		& 2/15 & 2.64E-02 & 2.04 & 2.87E-02 & 2.02  \\
		& 2/20 & 1.48E-02 & 2.01 & 1.61E-02 & 2.01  \\
		& 2/25 & 9.47E-03 & 2.01 & 1.03E-02 & 2.01 \\ [1ex] \hline		
		
		\multirow{4}{*}{PAVF} & 2/10  & 6.03E-02 & - & 6.49E-02 &  - \\
		& 2/15  & 2.63E-02 & 2.05 & 2.85E-02 & 2.02 \\
		& 2/20 & 1.48E-02 & 2.00 & 1.60E-02 & 2.02  \\
		& 2/25 & 9.42E-03 & 2.01 & 1.02E-02 & 2.01\\ [1ex]\hline	
		
		\multirow{4}{*}{PAVF-C} & 2/10  & 6.04E-02 & - & 6.50E-02 & - \\
		& 2/15  & 2.64E-02 & 2.04 & 2.87E-02 & 2.02 \\
		& 2/20 & 1.48E-02 & 2.01 & 1.61E-02 & 2.01 \\
		& 2/25 & 9.47E-03 & 2.01 & 1.03E-02 & 2.01 \\ [1ex]\hline	
		
		\multirow{4}{*}{PAVF-P} & 2/10  & 6.04E-02 & - & 6.50E-02 &  - \\
		& 2/15  & 2.64E-02 & 2.04 & 2.87E-02 & 2.02 \\
		& 2/20 & 1.48E-02 & 2.01 & 1.61E-02 & 2.01  \\
		& 2/25 & 9.47E-03 & 2.01 & 1.03E-02 & 2.01\\
		\hline
	\end{tabular*}
\end{table}

\begin{figure}[H]
	\centering
	(a)\includegraphics[width=0.45\textwidth]{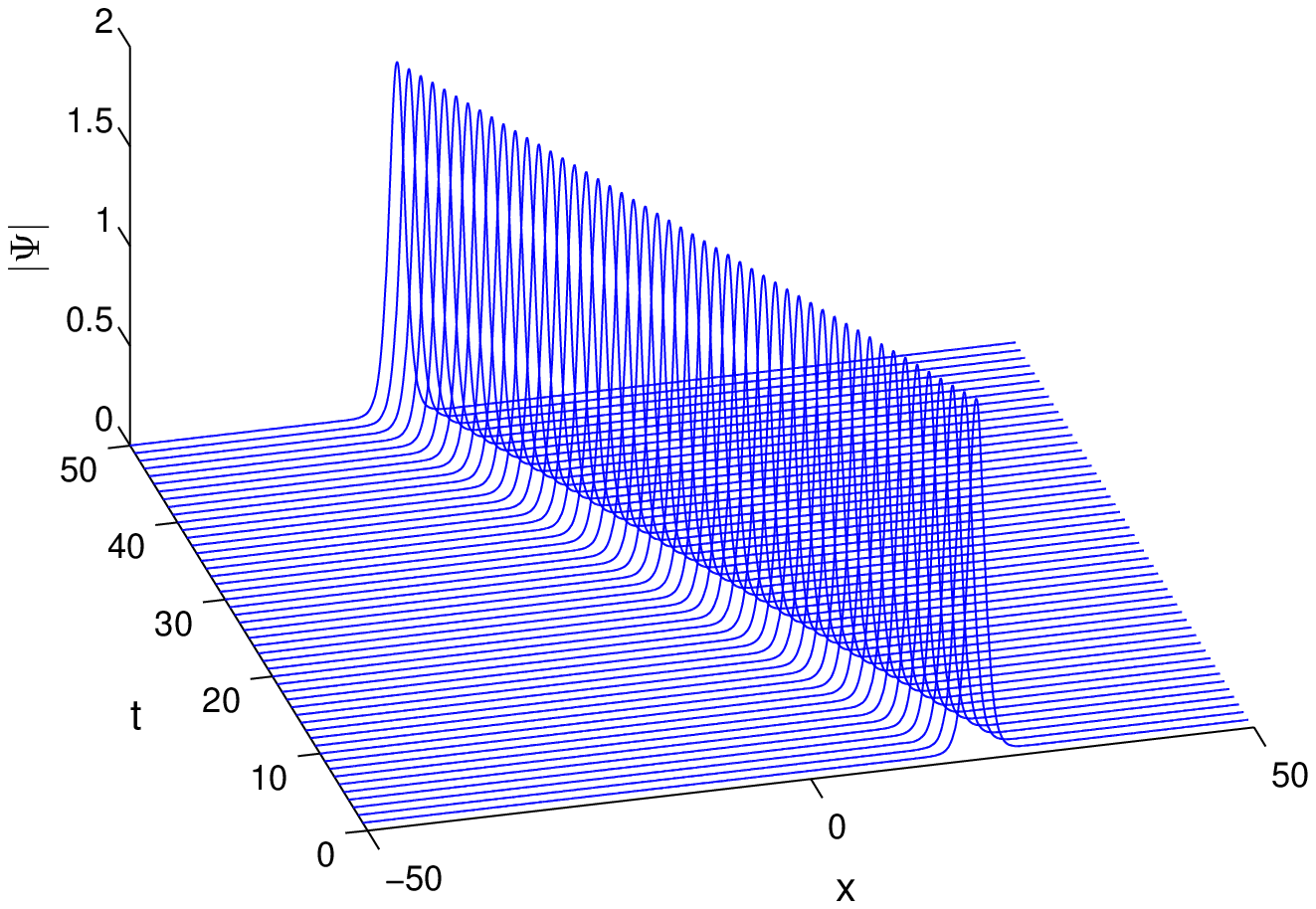}
	(b)\includegraphics[width=0.45\textwidth]{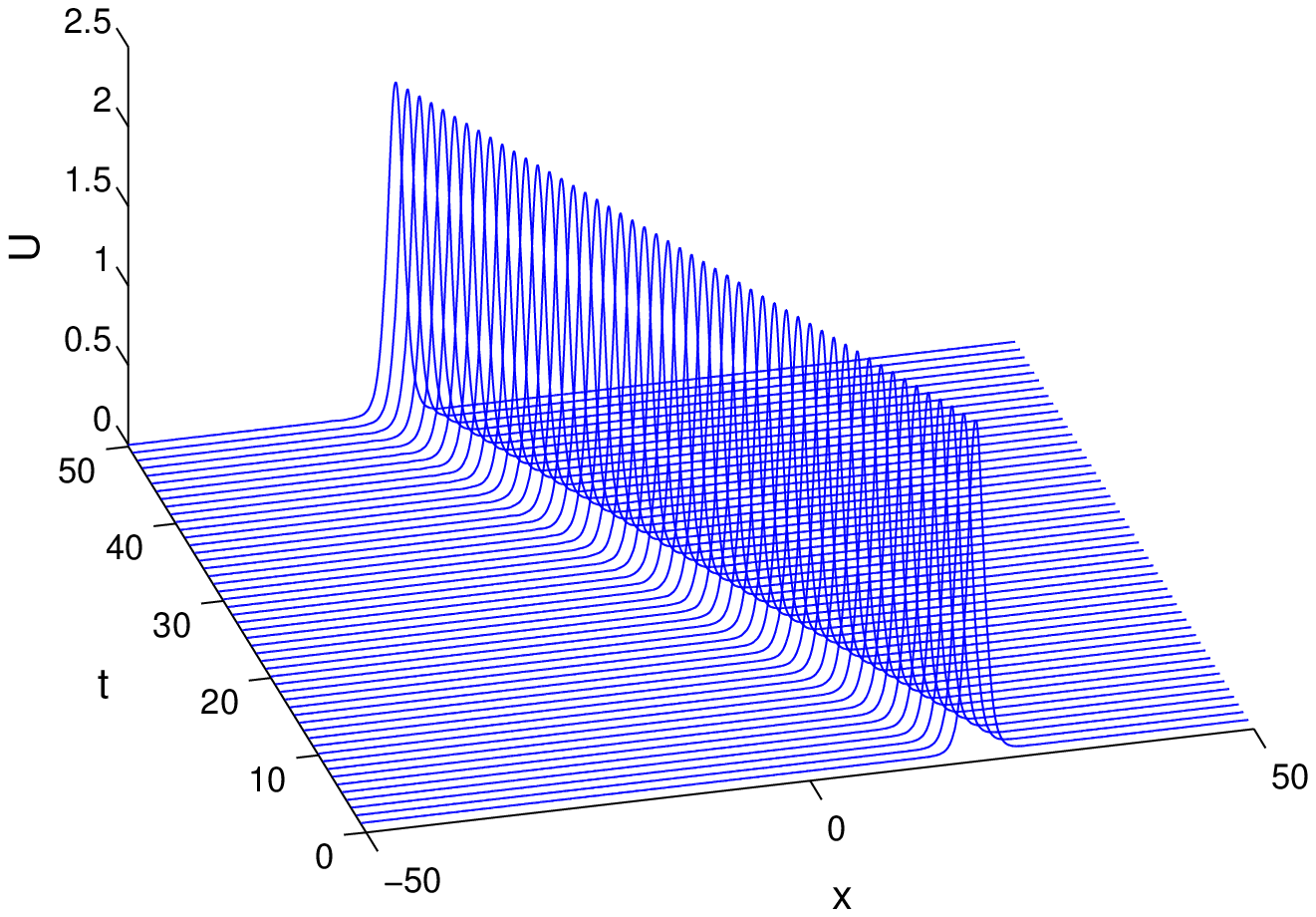}\\
	\caption{Profiles of one soliton generated by four energy-preserving methods.
		(a): $|Q+Pi|$. (b) $U$.}
	\label{fig:3-3}
\end{figure}

\begin{figure}[H]
	\centering
	(a)\includegraphics[width=0.45\textwidth]{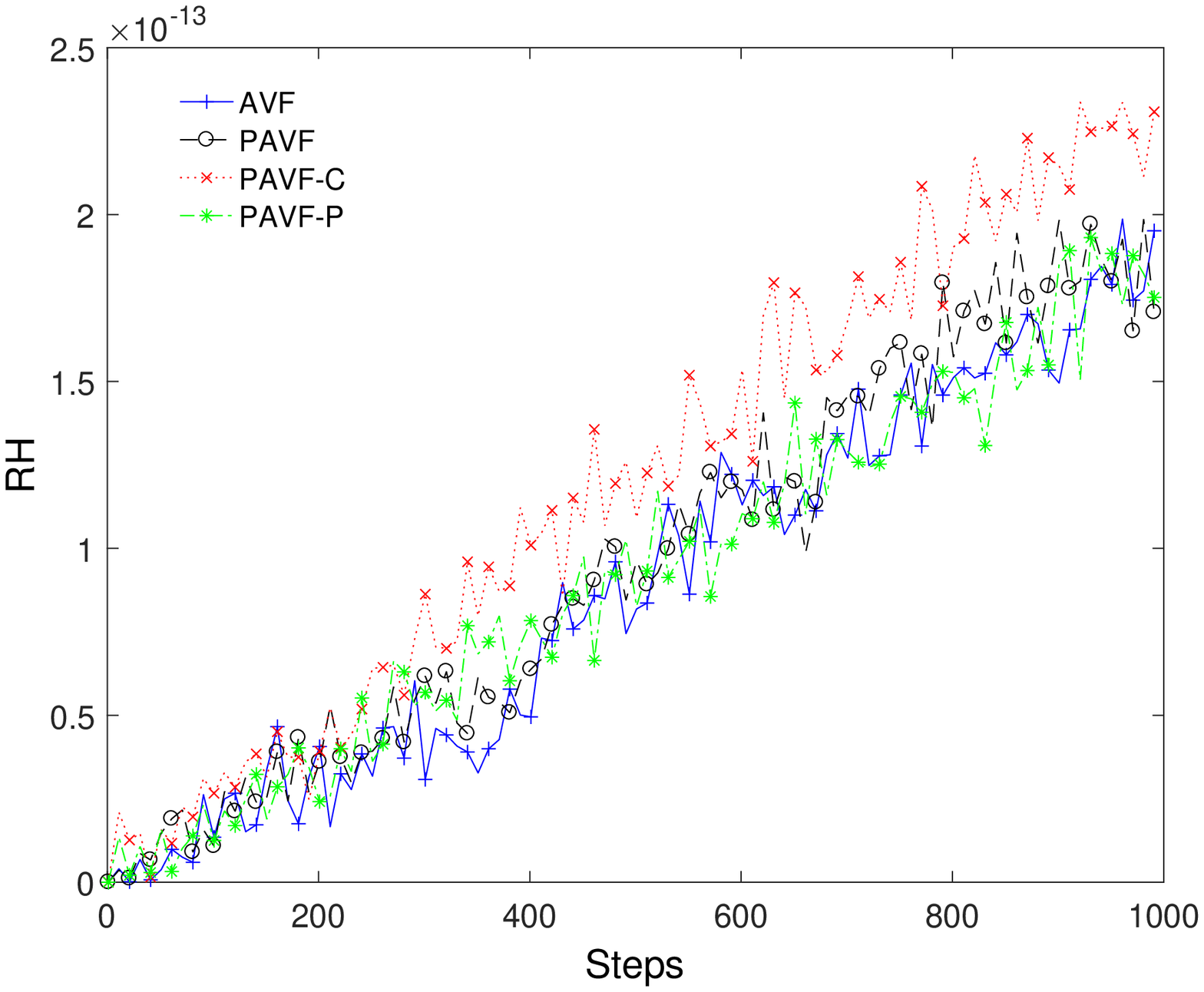}
	(b)\includegraphics[width=0.45\textwidth]{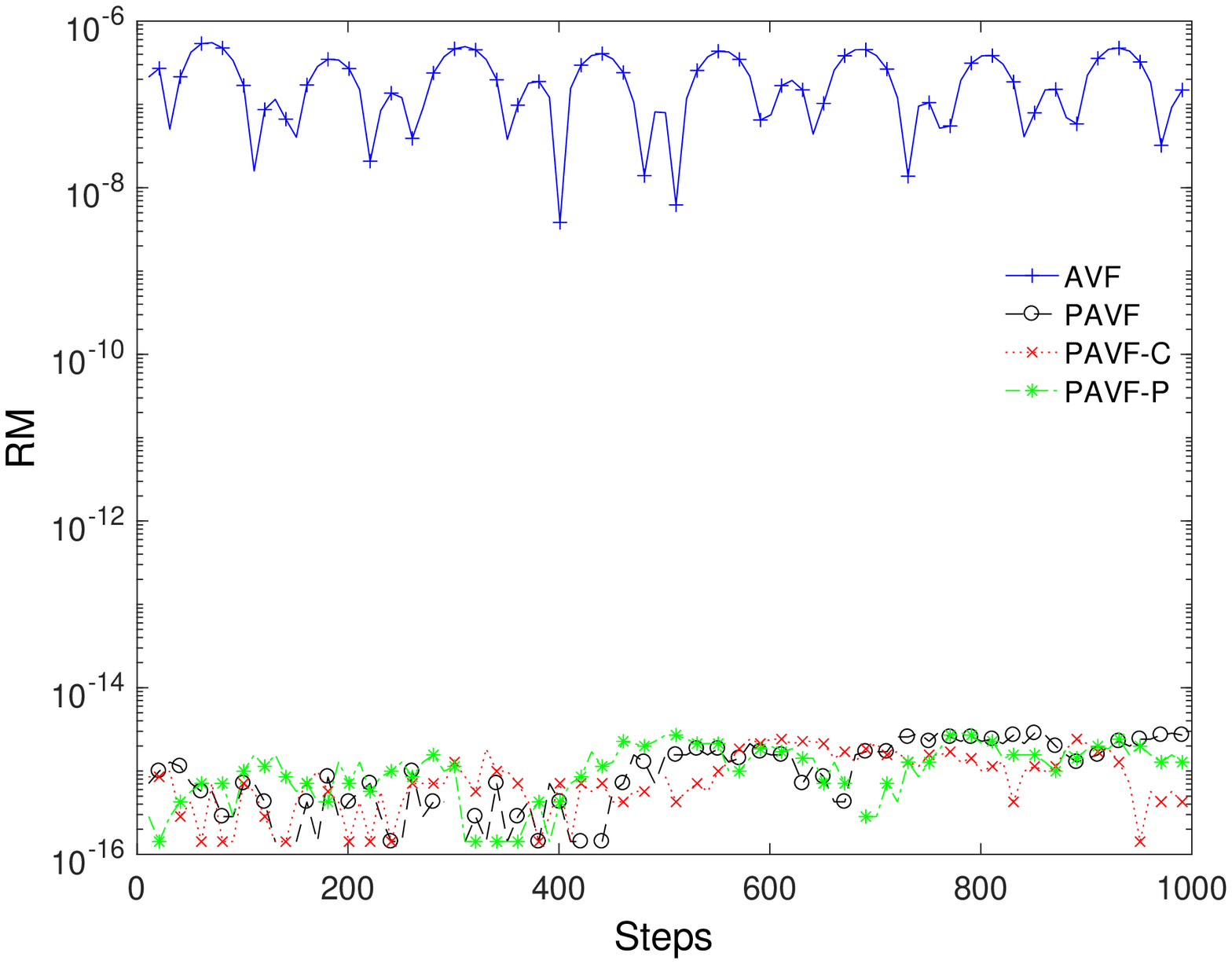}\\
	\caption{Relative energy error $RH$ (a) and mass error $RM$ (b) of  four energy-preserving methods for one soliton evolution.}
	\label{fig:3-4}
\end{figure}

To demonstrate the long-term behavior of the proposed methods, we enlarge the computational domain to $[-50,50]$. The spatial and temporal steps are set to $h=0.1$, $\tau=0.05$. The   initial phase in located at $x_{0}=20$. Figure.~\ref{fig:3-3} presents the wave profiles of $|\Psi|$ and $U$ from $t=0$ to $t=50$ which can be generated by all the four methods with almost same profiles. Moreover, in view of the relative errors in the energy and mass conservation laws (Figure.~\ref{fig:3-4}), we can find that the discrete energy can be preserved to round-off errors by any of the methods. However, the original AVF method can no long preserve the discrete mass conservation law even though the errors are always bounded during the time evolution. In contrast, the rest three partitioned AVF methods give an exact preservation of the discrete mass. In Table.~\ref{tab:3-4}, we list the CPU time of the four methods. Since the AVF method and the PAVF-P method are fully implicit, their computational costs are comparable but far more than the PAVF method and the PAVF-C method which are only linearly implicit. Moreover, the CPU time of the the PAVF-C method is much less than the twice of that of the PAVF method.

\begin{table}[H]
	\centering
	\caption{Computational cost of one solitons by four energy-preserving methods.}
	\label{tab:3-4}
	\begin{tabular*}{\textwidth}{@{\extracolsep{\fill}}ccccc}\hline
		$t=50$     & AVF                      & PAVF                      & PAVF-C                     & PAVF-P \\   \hline\noalign{\smallskip}
		CPU time (s)     & 183.44                  &   60.28                &     84.16        &    186.11  \\\hline
	\end{tabular*}
\end{table}

\begin{figure}[H]
	\centering
	(a)\includegraphics[width=0.45\textwidth]{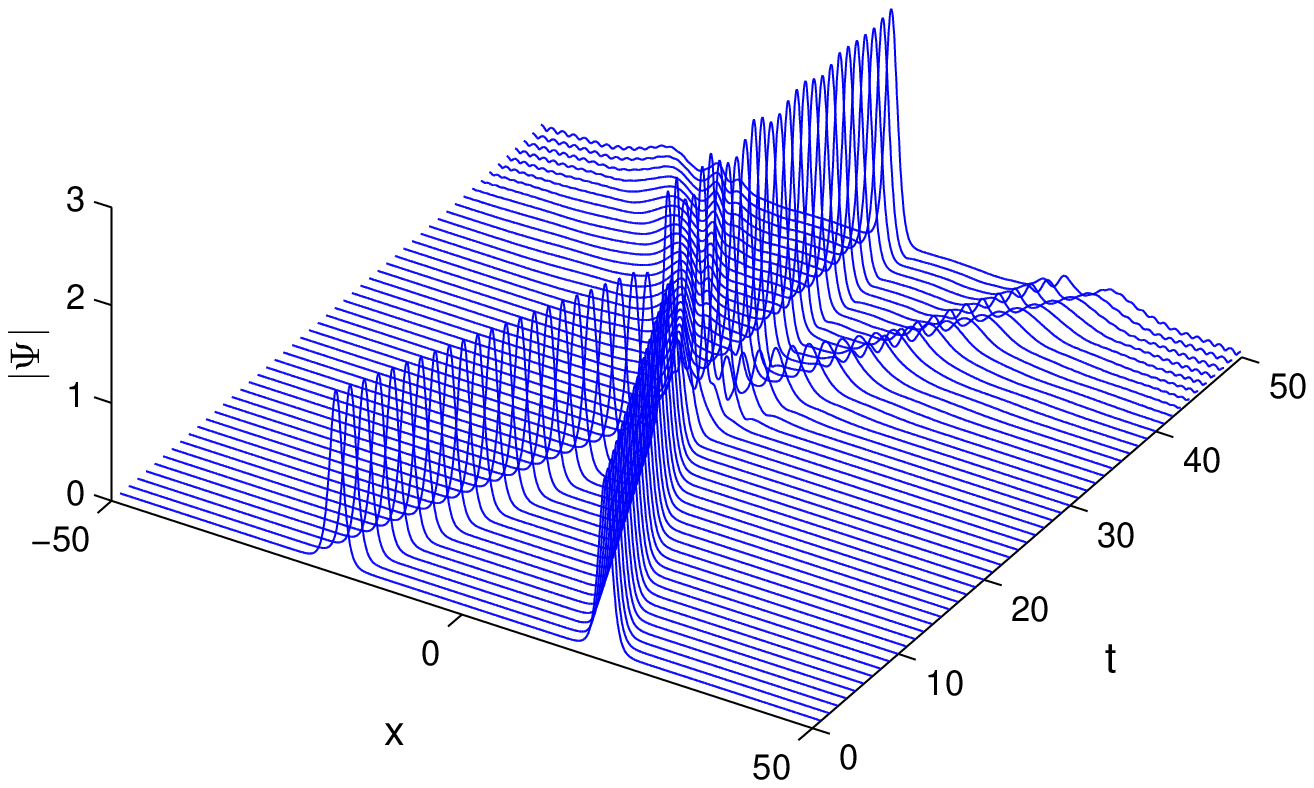}
	(b)\includegraphics[width=0.45\textwidth]{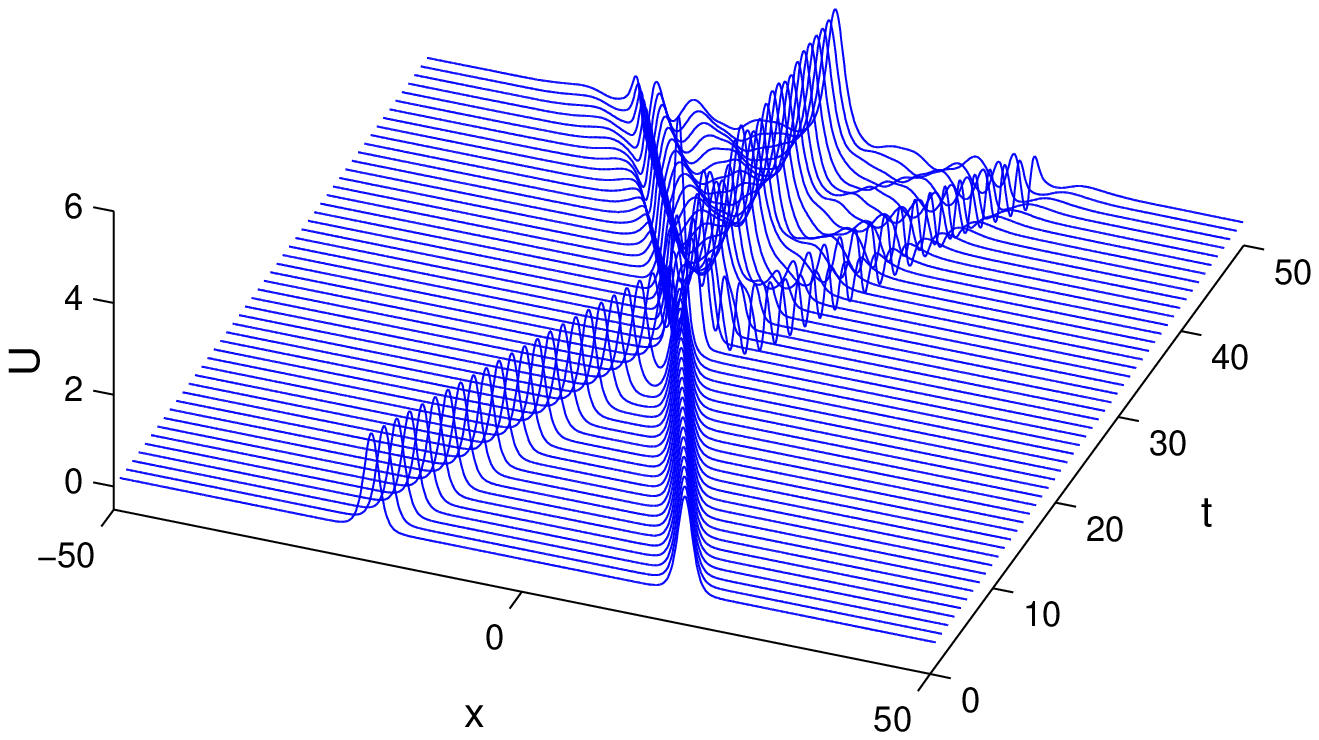}\\
	\caption{Collisions of two solitons generated by four energy-preserving methods.
		(a): $|\Psi|$. (b) $U$.}
	\label{fig:3-5}
\end{figure}

\begin{figure}[H]
	\centering
	(a)\includegraphics[width=0.45\textwidth]{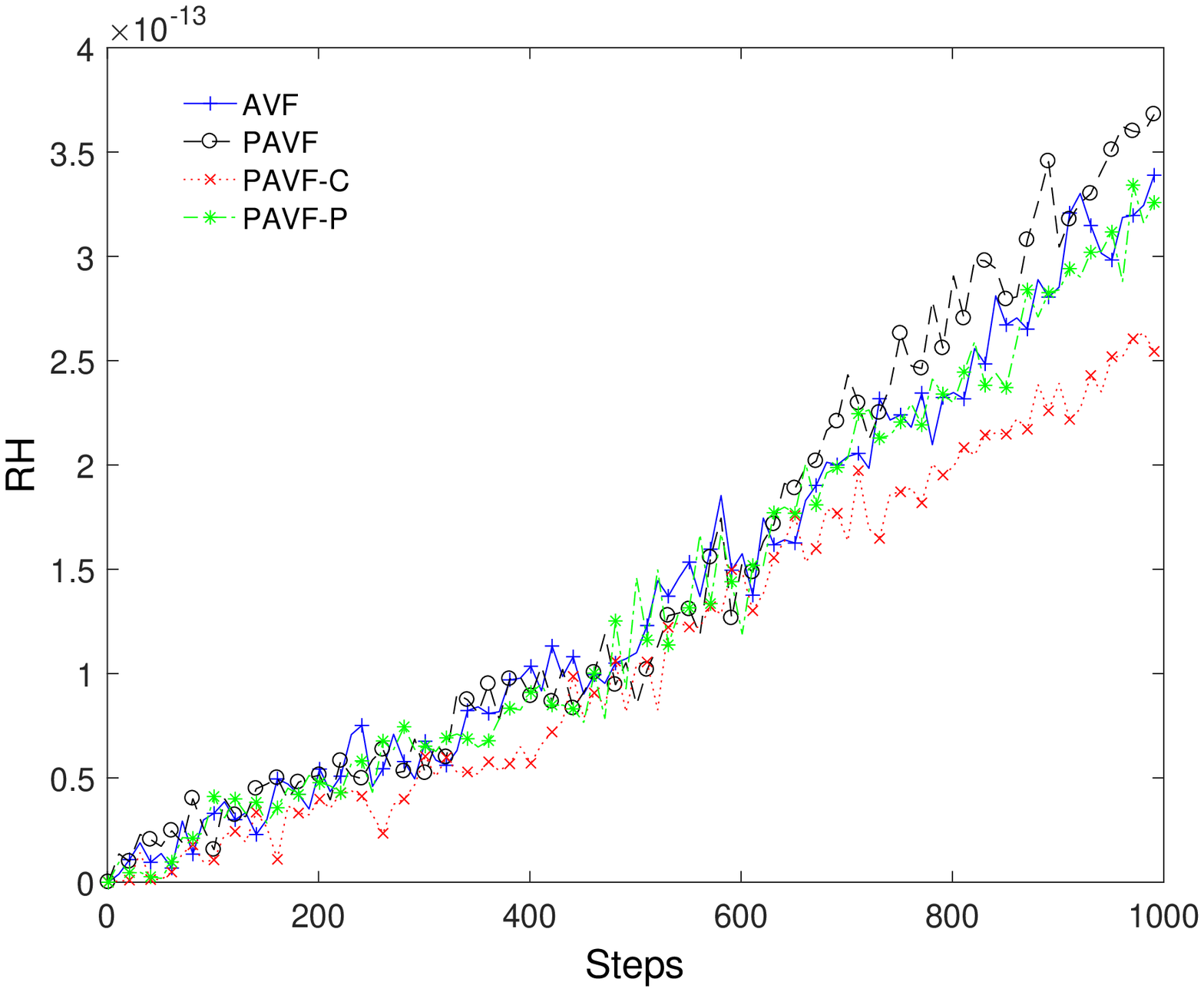}
	(b)\includegraphics[width=0.45\textwidth]{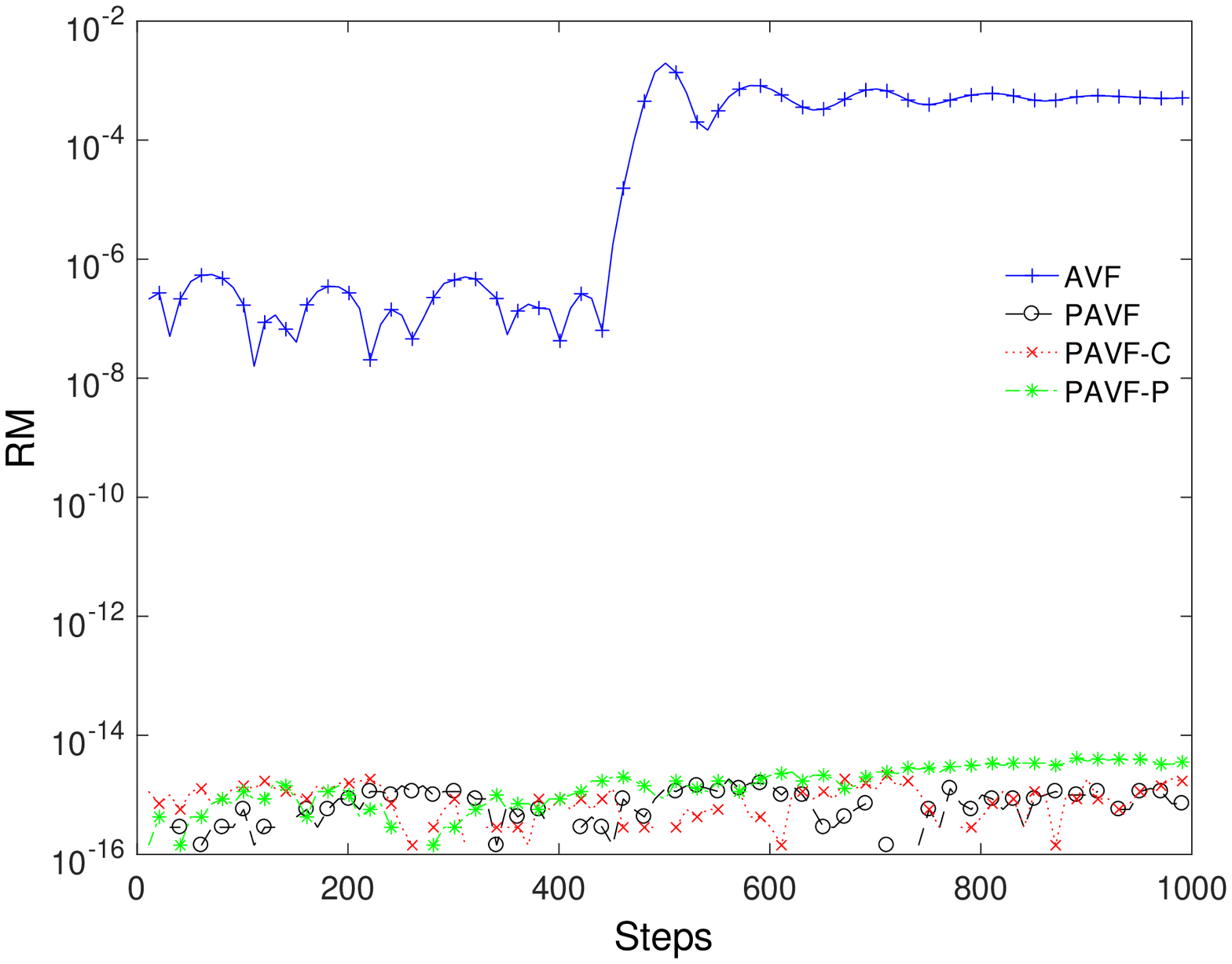}\\
	\caption{Relative energy error $RH$ (a) and mass error $RM$ (b) of the four energy-preserving methods for two solitons evolution.}
	\label{fig:3-6}
\end{figure}

Next we consider the case of two colliding solitons with initial condition
\begin{align*}
	&\varphi(x,0)=\sum_{i=1}^{2}\left(\frac{3\sqrt{2}}{4\sqrt{1-c_{i}^{2}}}\mathrm{sech}^{2}\Big(\frac{1}{2\sqrt{1-c_{i}^{2}}}(x-x_{i})\Big)\exp(ic_ix)\right),\\
	&u(x,0)=\sum_{i=1}^{2}\left(\frac{3}{4(1-c_{i}^{2})}\mathrm{sech}^{2}\Big(\frac{1}{2\sqrt{1-c_{i}^{2}}}(x-x_{i})\Big)\right),\\
	&u_{t}(x,0)=\sum_{i=1}^{2}\left(\frac{3c_{i}}{4(1-c_{i}^{2})^{3/2}}\mathrm{sech}^{2}\Big(\frac{1}{2\sqrt{1-c_{i}^{2}}}(x-x_{i})\Big)\tanh\Big(\frac{1}{2\sqrt{1-c_{i}^{2}}}(x-x_{i})\Big)\right),
\end{align*}
where $|c_i|<1$, $x_{i}$, $i=1,2$ are the propagating velocities and initial phases of two solitons, respectively. We set $c_1=-0.8, x_1=20$ and $c_2=0.8, x_2=-20$ which corresponds to two colliding solitons with same amplitudes and speed but opposite directions and initial phases. The spatial and temporal steps are $h=0.1, \tau=0.05$. Figure.~\ref{fig:3-5} demonstrates the evolution of shapes of $|\Psi|$ and $U$ until $t=50$ with clear collision observed. Although these profiles can be produced both by the original AVF method and our partitioned AVF methods, from the errors of discrete energy and mass in Figure.~\ref{fig:3-6} our methods are superior in the exact conservation of these two invariants while the original method can only preserve the discrete energy. The computational costs of four methods are listed in Table.~\ref{tab:3-5} which reveal a similar result as that in Table.~\ref{tab:3-4}.

\begin{table}[H]
	\centering
	\caption{Computational cost of two solitons by four energy-preserving methods.}
	\label{tab:3-5}
	\begin{tabular*}{\textwidth}{@{\extracolsep{\fill}}ccccc}\hline
		$t=50$     & AVF                      & PAVF                      & PAVF-C                     & PAVF-P \\   \hline\noalign{\smallskip}
		CPU time (s)     & 189.24                 &   57.80           &     84.18         & 185.80     \\\hline
	\end{tabular*}
\end{table}

\section{Conclusions}
\label{sec:4}

In this paper we present a partitioned AVF method for Hamiltonian ODEs and PDEs which differs from the conventional AVF method in the derivation of the vector of means of the tangential components in the gradient of the Hamiltonian. With the grouping strategy, the resulting schemes are semi-implicit or linearly implicit. Consequently, such schemes can decrease the computational scale or even avoid the iteration process which has to be implemented by the original AVF method. In addition, we further find that the partitioned AVF method can preserve extra conservative quantities except for the Hamiltonian energy in particular problem while the AVF method cannot. When considering all variables as one group, the PAVF method just becomes the standard AVF method. In another extreme when each variable is viewed as an individual group, then the PAVF method is actually equivalent to the discrete gradient method.

For sake of improving the accuracy of the original partitioned AVF method, in conjunction with its adjoint, we further introduce the partitioned AVF composition method and plus method. Both the two modifications inherit the energy conservative property. However, the partitioned AVF composition method is more efficient than the plus method. Therefore, in practical the partitioned AVF composition method is a perfect alternative to the conventional AVF method.

\vspace{0.3cm}
\hspace{-0.5cm}{\bf Acknowledgements}\\
This work is supported by the Jiangsu Collaborative Innovation Center for Climate
Change, the National Natural Science Foundation of China (Grant Nos. 11271195,
41231173, 41504078) and the Priority Academic Program Development of Jiangsu Higher Education Institutions.

\vspace{0.5cm}
\hspace{-0.5cm}{\bf \large Reference}

\end{document}